\documentclass[draftclsnofoot,journal,onecolumn, 12pt]{IEEEtran}
\usepackage{mathtools}
\usepackage{tabularx}
\usepackage{tikz}
\usetikzlibrary{arrows,positioning} 
\tikzset{
	to/.style={->,>=stealth',shorten >=1pt,semithick,font=\sffamily\footnotesize},
	every node/.style={align=center}
}

\usepackage[T1]{fontenc}
\usepackage{amsmath,amssymb,amsfonts,mathrsfs,bm}
\usepackage{mathtools}
\usepackage{amsthm}
\usepackage{cite}
\usepackage[shortlabels]{enumitem}
\usepackage{graphicx}
\usepackage{epstopdf}
\usepackage{url}
\usepackage{colortbl}
\usepackage{multirow}
\usepackage{xcolor}
\usepackage[normalem]{ulem}
\usepackage{array}
\newcolumntype{L}[1]{>{\raggedright\let\newline\\\arraybackslash\hspace{0pt}}m{#1}}
\newcolumntype{C}[1]{>{\centering\let\newline\\\arraybackslash\hspace{0pt}}m{#1}}
\newcolumntype{R}[1]{>{\raggedleft\let\newline\\\arraybackslash\hspace{0pt}}m{#1}}

\makeatletter
\let\MYcaption\@makecaption
\makeatother
\usepackage[font=footnotesize]{subcaption}
\makeatletter
\let\@makecaption\MYcaption
\makeatother

\usepackage{xparse}



\usepackage{glossaries}
\newacronym{wrt}{w.r.t.}{with respect to}
\newacronym{RHS}{R.H.S.}{right-hand side}
\newacronym{LHS}{L.H.S.}{left-hand side}
\newacronym{iid}{i.i.d.}{independent and identically distributed}

\usepackage{float}

\ifx\notloadhyperref\undefined
	\ifx\loadbibentry\undefined
		\usepackage[hidelinks,hypertexnames=false]{hyperref} 
	\else
		\usepackage{bibentry}
		\makeatletter\let\saved@bibitem\@bibitem\makeatother
		\usepackage[hidelinks,hypertexnames=false]{hyperref}
		\makeatletter\let\@bibitem\saved@bibitem\makeatother
	\fi
\else
	\relax
\fi

\usepackage[capitalize]{cleveref}
\crefname{equation}{}{}
\Crefname{equation}{}{}
\crefname{claim}{claim}{claims}
\crefname{step}{step}{steps}
\crefname{line}{line}{lines}
\crefname{dmath}{}{}
\crefname{dseries}{}{}
\crefname{dgroup}{}{}
\crefname{Theorem}{Theorem}{Theorems}
\crefname{Corollary}{Corollary}{Corollaries}
\crefname{Proposition}{Proposition}{Propositions}
\crefname{Lemma}{Lemma}{Lemmas}
\crefname{Definition}{Definition}{Definitions}
\crefname{Example}{Example}{Examples}
\crefname{Assumption}{Assumption}{Assumptions}
\crefname{Remark}{Remark}{Remarks}
\crefname{Rem}{Remark}{Remarks}
\crefname{remarks}{Remarks}{Remarks}
\crefname{Theorem_A}{Theorem}{Theorems}
\crefname{Corollary_A}{Corollary}{Corollaries}
\crefname{Proposition_A}{Proposition}{Propositions}
\crefname{Lemma_A}{Lemma}{Lemmas}
\crefname{Definition_A}{Definition}{Definitions}

\usepackage{algorithm,algorithmic}

\ifx\loadbreqn\undefined
	\relax
\else
	\usepackage{breqn} 
\fi

\ifx\renewtheorem\undefined
\ifx\useTheoremCounter\undefined
\newtheorem{Theorem}{Theorem}
\newtheorem{Corollary}{Corollary}
\newtheorem{Proposition}{Proposition}
\newtheorem{Lemma}{Lemma}
\else
\newtheorem{Theorem}{Theorem}

\newtheorem{Proposition}[theorem]{Proposition}
\fi

\newtheorem{Assumption}{Assumption}

\newtheorem{Lemma_A}{Lemma}[section]

\fi

\theoremstyle{remark}

\theoremstyle{plain}

\newcommand{\Real}{\mathbb{R}}
\newcommand{\Nat}{\mathbb{N}}

\newcommand{\calN}{\mathcal{N}}

\DeclareSymbolFont{bsfletters}{OT1}{cmss}{bx}{n}
\DeclareSymbolFont{ssfletters}{OT1}{cmss}{m}{n}
\DeclareMathSymbol{\bsfGamma}{0}{bsfletters}{'000}
\DeclareMathSymbol{\ssfGamma}{0}{ssfletters}{'000}
\DeclareMathSymbol{\bsfDelta}{0}{bsfletters}{'001}
\DeclareMathSymbol{\ssfDelta}{0}{ssfletters}{'001}
\DeclareMathSymbol{\bsfTheta}{0}{bsfletters}{'002}
\DeclareMathSymbol{\ssfTheta}{0}{ssfletters}{'002}
\DeclareMathSymbol{\bsfLambda}{0}{bsfletters}{'003}
\DeclareMathSymbol{\ssfLambda}{0}{ssfletters}{'003}
\DeclareMathSymbol{\bsfXi}{0}{bsfletters}{'004}
\DeclareMathSymbol{\ssfXi}{0}{ssfletters}{'004}
\DeclareMathSymbol{\bsfPi}{0}{bsfletters}{'005}
\DeclareMathSymbol{\ssfPi}{0}{ssfletters}{'005}
\DeclareMathSymbol{\bsfSigma}{0}{bsfletters}{'006}
\DeclareMathSymbol{\ssfSigma}{0}{ssfletters}{'006}
\DeclareMathSymbol{\bsfUpsilon}{0}{bsfletters}{'007}
\DeclareMathSymbol{\ssfUpsilon}{0}{ssfletters}{'007}
\DeclareMathSymbol{\bsfPhi}{0}{bsfletters}{'010}
\DeclareMathSymbol{\ssfPhi}{0}{ssfletters}{'010}
\DeclareMathSymbol{\bsfPsi}{0}{bsfletters}{'011}
\DeclareMathSymbol{\ssfPsi}{0}{ssfletters}{'011}
\DeclareMathSymbol{\bsfOmega}{0}{bsfletters}{'012}
\DeclareMathSymbol{\ssfOmega}{0}{ssfletters}{'012}

\DeclareMathOperator*{\argmax}{arg\,max}

\DeclareMathOperator{\st}{s.t.}

\DeclareMathOperator*{\esssup}{ess\,sup}

\newcommand{\qednew}{\nobreak \ifvmode \relax \else
      \ifdim\lastskip<1.5em \hskip-\lastskip
      \hskip1.5em plus0em minus0.5em \fi \nobreak
      \vrule height0.75em width0.5em depth0.25em\fi}

\newcommand{\nn}{\nonumber\\}

\newcommand{\ud}{\mathrm{d}}

\newcommand{\indicator}[1]{{\bf 1}_{\{{#1}\}}}
\newcommand{\indicatore}[1]{{\bf 1}_{#1}}

\newcommand{\ofrac}[1]{{\frac{1}{#1}}}

\newcommand{\ceil}[1]{\left\lceil{#1}\right\rceil}
\newcommand{\floor}[1]{\left\lfloor{#1}\right\rfloor}

\newcommand{\KLD}[2]{{D({#1}\, ||\, {#2})}}

\newcommand{\cond}[2]{\left. {#1}\, \middle| \, {#2} \right.}

\DeclareDocumentCommand \P { g d() g } {%
	\IfNoValueTF {#3} 
	{%
		\IfNoValueTF {#1} 
		{%
			\IfNoValueTF {#2}
			{%
				\mathbb{P}%
			}%
			{%
				\mathbb{P}\left({#2}\right)%
			}%
		}%
		{%
			\IfNoValueTF {#2}
			{%
				\mathbb{P}_{#1}%
			}%
			{%
				\mathbb{P}_{#1}\left({#2}\right)%
			}%
		}%
	}%
	{%
		\IfNoValueTF {#1} 
		{%
			\mathbb{P}\left(\cond{#2}{#3}\right)%
		}%
		{%
			\mathbb{P}_{#1}\left(\cond{#2}{#3}\right)%
		}%
	}%
}

\DeclareDocumentCommand \E { g o g } {%
	\IfNoValueTF {#3} 
	{%
		\IfNoValueTF {#1} 
		{%
			\IfNoValueTF {#2}
			{%
				\mathbb{E}%
			}%
			{%
				\mathbb{E}\left[{#2}\right]%
			}%
		}%
		{%
			\IfNoValueTF {#2}
			{%
				\mathbb{E}_{#1}%
			}%
			{%
				\mathbb{E}_{#1}\left[{#2}\right]%
			}%
		}%
	}%
	{%
		\IfNoValueTF {#1} 
		{%
			\mathbb{E}\left[\cond{#2}{#3}\right]%
		}%
		{%
			\mathbb{E}_{#1}\left[\cond{#2}{#3}\right]%
		}%
	}%
}

\definecolor{gray90}{gray}{0.9}

\ifx\nohighlights\undefined

	\newcommand{\msout}[1]{\text{\color{green} \sout{\ensuremath{#1}}}}
	\newcommand{\del}[1]{{\color{green}\ifmmode \msout{#1}\else\sout{#1}\fi}}
\else

	\newcommand{\msout}[1]{#1}
	\newcommand{\del}[1]{#1}
\fi

\newcommand{\hide}[1]{}

\renewcommand{\figurename}{Fig.}
\newcommand{\figref}[1]{\figurename~\ref{#1}}
\graphicspath{{./Figures/}} 
\ifx\diagnoselabel\undefined
	\relax
\else
	\makeatletter
	 \def\@testdef #1#2#3{%
		 \def\reserved@a{#3}\expandafter \ifx \csname #1@#2\endcsname
		\reserved@a  \else
	 \typeout{^^Jlabel #2 changed:^^J%
	 \meaning\reserved@a^^J%
	 \expandafter\meaning\csname #1@#2\endcsname^^J}%
	 \@tempswatrue \fi}
	\makeatother
\fi

\newacronym{GLR}{GLR}{generalized likelihood-ratio}
\newacronym{GLRT}{GLRT}{Generalized Likelihood Ratio Test}
\newacronym{ARL}{ARL}{average run length}
\newacronym{WADD}{WADD}{worst-case average detection delay}
\newacronym{ADD}{ADD}{average detection delay}
\newacronym{TCD}{TCD}{transient change detection}
\newacronym{KL}{KL}{Kullback-Leibler}
\newacronym{FMA}{FMA}{finite moving average}

\graphicspath{{./Figures/}}

\newcommand{\ARL}{\text{ARL}}
\newcommand{\WADD}{\text{WADD}}

\newcommand{\tauWSGLR}{\tau_{\text{W-SGLR}}}
\newcommand{\tauSGLR}{\tau_{\text{SGLR}}}

\title{Quickest Change Detection in the Presence of a Nuisance Change}

\author{
	Tze~Siong~Lau, and Wee~Peng~Tay,~\IEEEmembership{Senior Member,~IEEE}%
	\thanks{This research is supported by the Singapore Ministry of Education Academic Research Fund Tier 1 grant 2017-T1-001-059 (RG20/17) and Tier 2 grant MOE2018-T2-2-019.}%
	\thanks{T.~S.~Lau and W.~P.~Tay are with the School of Electrical and Electronic Engineering, Nanyang Technological University, Singapore (e-mail: TLAU001@e.ntu.edu.sg, wptay@ntu.edu.sg). 
	}
}

\begin{document}
	\maketitle
	\begin{abstract}
		In the quickest change detection problem in which both nuisance and critical changes may occur, the objective is to detect the critical change as quickly as possible without raising an alarm when either there is no change or a nuisance change has occurred. A window-limited sequential change detection procedure based on the generalized likelihood ratio test statistic is proposed. A recursive update scheme for the proposed test statistic is developed and is shown to be asymptotically optimal under mild technical conditions. In the scenario where the post-change distribution belongs to a parametrized family, a generalized stopping time and a lower bound on its average run length are derived. The proposed stopping rule is compared with the \gls{FMA} stopping time and the naive 2-stage procedure that detects the nuisance or critical change using separate CuSum stopping procedures for the nuisance and critical changes. Simulations demonstrate that the proposed rule outperforms the \gls{FMA} stopping time and the 2-stage procedure, and experiments on a real dataset on bearing failure verify the performance of the proposed stopping time.
	\end{abstract}

	\section{Introduction}\label{sec:intro}
	
	The problem of detecting a change in the statistical properties of a signal with the shortest possible delay after the change is known as quickest change detection (QCD). Given a sequence of independent and identically distributed (i.i.d.) observations $\{x_t:t\in \mathbb{N}\}$ with distribution $f$ up to an unknown change point $\nu$ and i.i.d.\ with distribution $g\neq f $ after $\nu$, we aim to detect this change as quickly as possible while maintaining a false alarm constraint. Detecting for a change has applications in many areas, including manufacturing quality control\cite{woodall2004using,lai1995sequential}, fraud detection\cite{bolton02}, cognitive radio\cite{lai2008quickest}, network surveillance\cite{sequeira02,tartakovsky2006novel,LuoTayLen:J16}, structural health monitoring\cite{sohn00}, spam detection\cite{xie12,JiTayVar:J17,TanJiTay:J18}, bioinformatics\cite{muggeo10}, power system line outage detection\cite{banerjee2014power}, and sensor networks\cite{coppin1996digital,Hong2004,YanZhoTay:J18}.
	
	
	For the non-Bayesian formulation of QCD, the change-point is assumed to be unknown but deterministic. When both the pre- and post-change distributions are known, Page \cite{page54} developed the Cumulative Sum Control Chart (CuSum) for quickest change detection. Lorden\cite{lorden71} proved that the CuSum test has asymptotically optimal worst-case average detection delay as the false alarm rate goes to zero. Moustakides \cite{moustakides86} later established that the CuSum test is exactly optimal under Lorden's optimality criterion. Later,  Lai showed in \cite{lai98} that the CuSum test is asymptotically optimal under Pollak's criterion\cite{pollak85}, as the false alarm rate goes to zero. For the case where the post-change distribution is unknown, Lorden\cite{lorden71} showed that the \gls{GLR} CuSum is asymptotically optimal for the case of finite multiple post-change distributions. Other methods were also proposed for the case when the post-change distribution is unknown to a certain degree \cite{siegmund95,lai98,banerjee15,lau17,lau2017optimal,lau19}. We refer the reader to \cite{tartakovsky2014sequential,veeravalli2013quickest,poor2009quickest} and the references therein for an overview of the QCD problem.
	
	In many practical applications, the signal of interest may undergo different types of change. However, only a subset of these changes may be of interest to the user. One example is the problem of bearing failure detection using accelerometer readings\cite{smith2015rolling}. During normal operations, the bearings are driven at two different activity levels, idle or active. In a typical bearing failure detection scenario, the bearing is initially driven at the idle state. A change to  the active state results in a change in the statistical properties of the accelerometer readings. However, this change is not of interest to us and is called a nuisance change. We are only interested in the change arising from the failure of the bearing, which is known as a critical change. Furthermore, the statistical properties of the observations obtained when the bearing is faulty depend on the activity level that it is driven at. The traditional QCD framework does not allow us to distinguish between critical and nuisance changes. Furthermore, due to the nuisance change, the observations are no longer i.i.d.\ either in the pre-change or post-change regime, depending on when the nuisance change occurs. In this paper, we investigate the non-Bayesian formulation of the QCD problem under a nuisance change, and propose a window-limited stopping time that ignores the nuisance change but detects the critical change as quickly as possible.
	
	\subsection{Related Work}
	
	Existing works in QCD that consider the problem where observations are not generated i.i.d.\ before and after the change-point can be categorized into three main categories. In the first category, the papers\cite{fuh2003,fuh2004asymptotic,fuh2015quickest} consider the problem where the pre-change distribution and the post-change distribution are modeled as hidden Markov models (HMMs). In \cite{fuh2003}, the authors proved the asymptotic optimality of the CuSum procedure for the HMM signal model in the sense of Lorden. In \cite{fuh2004asymptotic}, the authors developed the Shiryayev-Roberts-Pollak (SRP) rule for the HMM signal model and proved its optimality in the sense of Pollak. The authors of \cite{fuh2015quickest} consider the problem where the vector parameter of a two-state HMM changes at some unknown time. The second category of papers\cite{moustakides1998quickest,tartakovsky2005general} considers a QCD problem which relaxes the i.i.d.\ assumption. In \cite{moustakides1998quickest}, the authors established the optimality of CuSum and the Shiryayev-Roberts stopping rule in the class of random processes with likelihood ratios that satisfy certain independence and stationary conditions. The class of random processes includes Markov chains, AR processes, and processes evolving on a circle. In \cite{tartakovsky2005general}, the authors considered the Bayesian QCD problem where conditions on the asymptotic behavior of the likelihood process are assumed. Unlike all the aforementioned papers, the signal model in our QCD problem with nuisance change cannot be modeled by an HMM, and the likelihood ratios generated by our signal model are non-stationary. In the third category, the papers \cite{Krishnamurthy12,guepie2012sequential,ebrahimzadeh2015sequential,Moustakides16,zou2018quickest,heydari2018quickest} consider QCD of transient changes, where the change is either not persistent or multiple changes occur throughout the monitoring process. Unlike our QCD problem which allows some changes to be considered nuisance, all changes are considered critical in the aforementioned papers.
	
	\subsection{Our Contributions}\label{subsec:contributions}
	
	In this paper, we consider the non-Bayesian QCD problem where both nuisance and critical changes may occur, and our objective is to detect the critical change as quickly as possible while ignoring the nuisance change. Our goal is to develop a sequential algorithm with computational complexity that increases linearly with the number of samples observed. Our main contributions are as follows:
	
	\begin{enumerate}
		\item We formulate the QCD problem with a nuisance change and propose a window-limited simplified \gls{GLR} (W-SGLR) stopping time.
		\item We derive a lower bound for the \gls{ARL} to a false alarm, and the asymptotic upper bound of the \gls{WADD} for our proposed test.
		\item We prove the asymptotic optimality of the W-SGLR stopping time under mild technical assumptions. 
		\item We provide simulation and experimental results that verify the theoretical guarantees of our proposed test and also illustrate the performance of our proposed test on a real dataset.
	\end{enumerate}
	
	A preliminary version of this work was presented in \cite{lau2018quickest,lau2019quickest}. To the best of our knowledge, there are no existing works that consider the QCD problem for a signal that may undergo a nuisance change.
	
	The rest of this paper is organized as follows. In Section~\ref{sec:problem}, we present our signal model and problem formulation. We propose the W-SGLR stopping time and derive the theoretical properties of our test statistics in Section~\ref{sec:teststats}. In Section~\ref{sec:discussion}, we discuss a modification of the proposed stopping time when the post-change distribution belongs to a parametrized family. We present numerical simulations and experiments on a real dataset to illustrate the performance of our proposed stopping time in Section~\ref{sec:numerical}. We conclude in Section~\ref{sec:conclusion}.

	\emph{Notations:} The operator $\E_f$ denotes mathematical expectation \gls{wrt} the probability density (pdf) $f$, and $X \sim f$ means that the random variable $X$ has distribution with pdf $f$. If the nuisance change point is at $\nu_n$, and the critical change point is at $\nu_c$, we let $\P{\nu_n,\nu_c}$ and $\E_{\nu_n,\nu_c}$ be the probability measure and mathematical expectation, respectively. The Gaussian distribution with mean $\mu$ and variance $\sigma^2$ is denoted as $\calN(\mu,\sigma^2)$. Convergence in $\P$-probability is denoted as $\xrightarrow{\ \P\ }$. We use $\indicatore{E}$ as the indicator function of the set $E$, and $\KLD{\cdot}{\cdot}$ to denote the \gls{KL} divergence. We use $\Nat$, $\mathbb{R}$ and $\mathbb{R}_{>0}$ to denote the set of positive integers, real numbers and positive real numbers, respectively.
	
	\section{Problem formulation}\label{sec:problem}
	
	In many applications, the statistical distribution of the observed signal may undergo different changes over time. For example, in the application of fault detection in motor bearings \cite{smith2015rolling}, we aim to raise an alarm as soon as possible after a bearing fault has occurred (critical change). This is done by monitoring the accelerometer readings from the motor to detect any changes in the signal statistics. However, the accelerometer readings are also affected by non-critical or nuisance changes like variation in the motor-load of the bearing. It would be undesirable if we declare that a fault has taken place whenever the motor-load changes. This motivates a need to define a change-point model that allows both critical and nuisances changes and to develop change detection techniques that can effectively ignore nuisance changes while efficiently identifying critical changes.

\begin{figure}
		\centering
		\begin{subfigure}{.25\textwidth}
			\centering
			\includegraphics[width=0.8\linewidth]{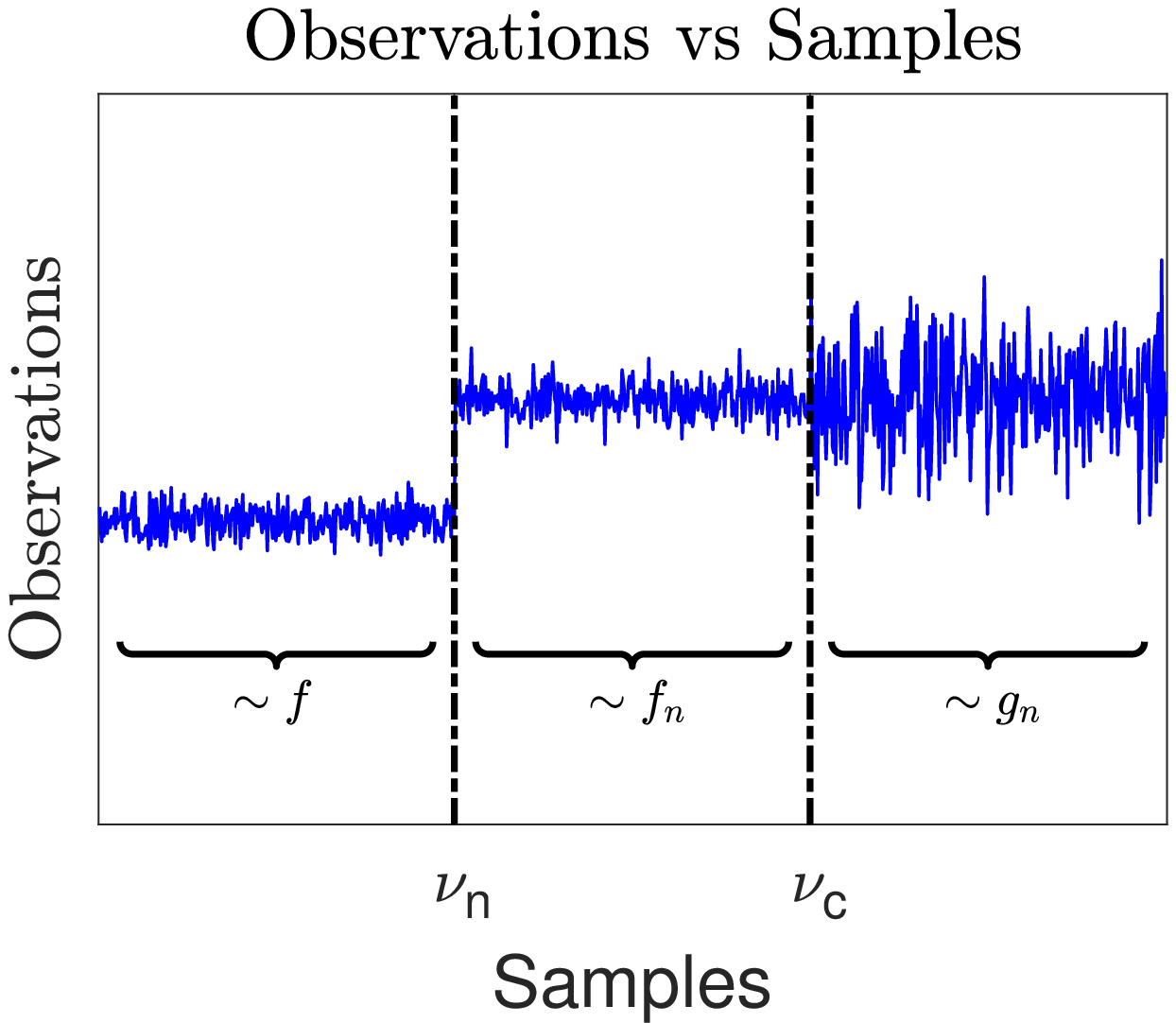}
			\caption{$\nu_n<\nu_c$.}
			\label{fig:sub1}
		\end{subfigure}%
		\begin{subfigure}{.25\textwidth}
			\centering
			\includegraphics[width=0.8\linewidth]{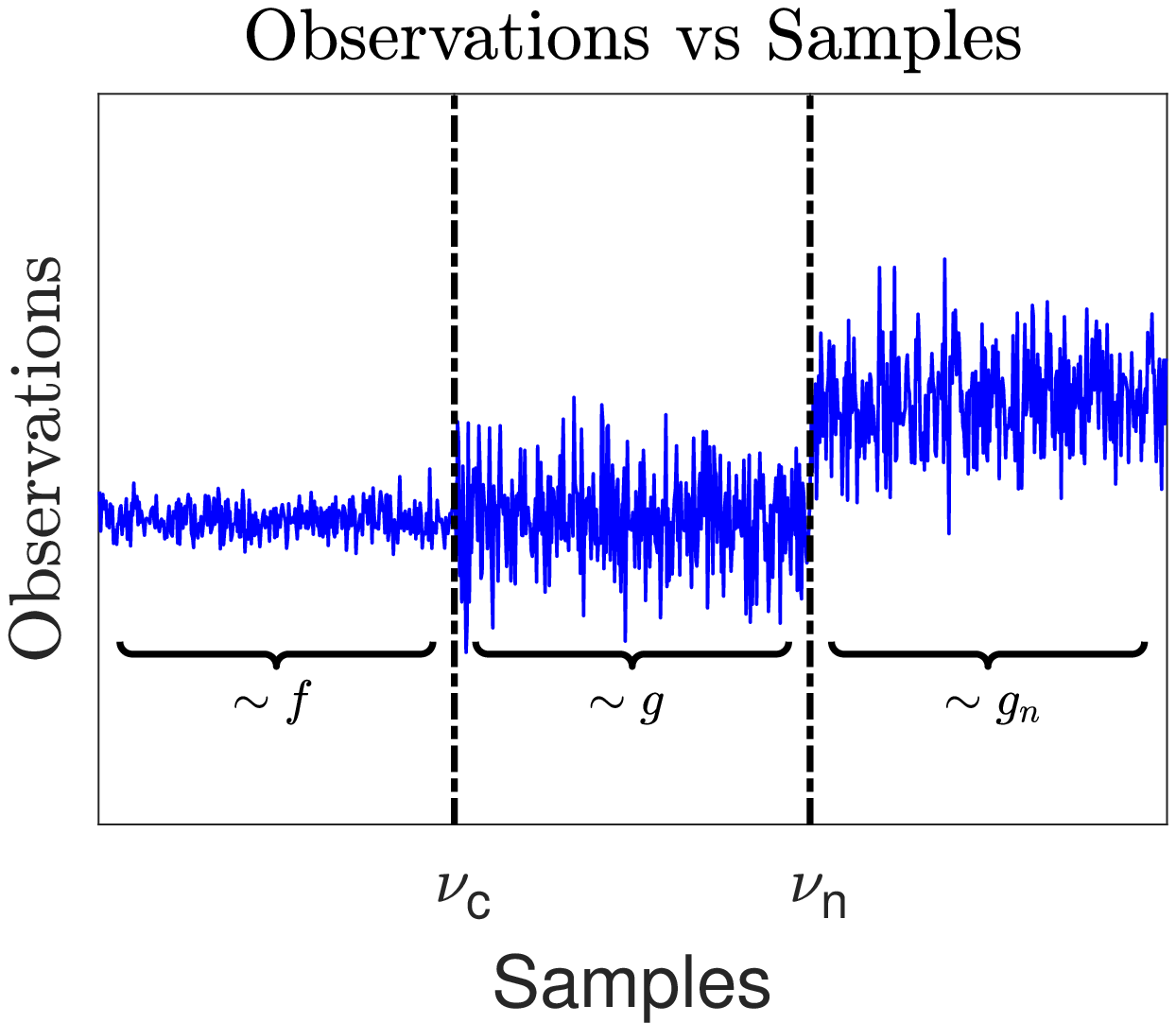}
			\caption{$\nu_n>\nu_c$.}
			\label{fig:sub2}
		\end{subfigure}
		\caption{Evolution of the signal distribution in our model.}	
		\label{fig:h_model}
	\end{figure}

	In this paper, we assume that the signals observed, $X_1, X_2, \ldots$, may undergo two types of change: a critical change at $\nu_c \geq 0$ and a nuisance change at $\nu_n \geq 0$. Both the critical and nuisance change points are unknown \emph{a priori}. We are interested in detecting the critical change while the nuisance change is not of interest. Let $f,f_n,g,g_n$ be probability distributions. At each time $t$, we let $h_{\nu_n,\nu_c,t}$ to be the distribution that generates the observation $X_t$ when the nuisance change point is at $\nu_n$ and the critical change point is at $\nu_c$:
	\begin{align}\label{eqn:definition_of_h}
	h_{\nu_n,\nu_c,t}=\begin{cases}
	f\quad \text{if $t<\min\{\nu_c,\nu_n\}$,}       \\
	f_n\quad \text{if $\nu_n\leq t<\nu_c$,}         \\
	g\quad \text{if $\nu_c\leq t<\nu_n$,}           \\
	g_n\quad \text{if $\max\{\nu_c,\nu_n\}\leq t$.} \\
	\end{cases}
	\end{align}
	Thus, in our model (cf. \cref{fig:h_model}), $f$ is the pre-change distribution.If $\nu_c < \nu_n < \infty$, the signal distribution first changes to $g$ at $\nu_c$ and then to $g_n$ at $\nu_n$. If $\nu_n < \nu_c < \infty$, the distribution first changes to $f_n$ at $\nu_n$ and then to $g_n$ at $\nu_c$. If $\nu_n=\nu_c$, then the distribution changes from $f$ to $g_n$ at the common change point.
	
	The sequence of observations $X_1,X_2,\ldots$ is a sequence of random variables satisfying $X_t \sim h_{\nu_n,\nu_c,t}$ where $\{X_t\}_{t\in\mathbb{N}}$ are mutually independent given $\nu_n,\nu_c$. 
	The quickest change detection problem is to detect the critical change $\nu_c$ through observing $X_1,X_2,\ldots,$ as quickly as possible while ignoring the nuisance change and keeping the false alarm rate low. In our signal model, the nuisance change also changes the distribution that generates the observations after the critical change point. This creates a dependence between the nuisance change point and the distribution after the critical change point. Our formulation is different from assuming composite pre-change and post-change distribution families\cite{mei2003asymptotically} since the nuisance change leads to non-stationarity in the distribution of $X_t$ before or after the critical change, depending on whether the nuisance change occurs before or after the critical change, respectively.
	
	In a typical sequential change detection procedure, at each time $t$, a test statistic $S(t)$ is
	computed based on the currently available observations $X_1,\ldots,X_t$, and the observer decides that a change has occurred at a stopping time
	$
	\tau(b)=\inf\{t:S(t)\geq b\}.
	$
	
	In the traditional QCD framework\cite{lorden71}, the rate of false alarms is quantified by the mean time between false alarms. Since the nuisance change-point affects the distributions generating the signal, this quantity varies with the nuisance change point. In this paper, we consider the worst-possible rate as the nuisance change point varies by considering the smallest mean time between false alarms for all possible nuisance change points. A similar generalization can be made to quantify the detection delay by taking the largest detection delay over all possible nuisance change points.  
	
	Mathematically, our QCD problem can be formulated as a minimax problem similar to Lorden's formulation\cite{lorden71}, where we seek a stopping time that minimizes the WADD subject to an \gls{ARL} constraint:
	\begin{equation}\label{eqn:lorden_formulation}
	\begin{aligned}
	& \min_\tau &  & \WADD(\tau)           \\
	& \st       &  & \ARL(\tau)\geq\gamma,
	\end{aligned}
	\end{equation}
	where $\gamma$ is a predefined threshold, $\tau$ is a stopping time \gls{wrt} the filtration $\{\sigma(X_1,X_2,\ldots,X_t) :\ t\geq 0\}$,\
	
	\begin{align}
	&\WADD(\tau) =\sup_{\nu_n\in\mathbf{N}\cup\{\infty\}}\WADD_{\nu_n}(\tau),\label{eqn:WADD}\\
	&\WADD_{\nu_n}(\tau) = \sup_{\nu_c\geq 1}\esssup \E{\nu_n,\nu_c}[(\tau-\nu_c+1)^+]{X_1,\ldots,X_{\nu_c-1}},\\
	&\ARL(\tau) =\inf_{\nu_n\in\mathbb{N}\cup\{\infty\}}\E{\nu_n,\infty}[\tau],
	\end{align}
	and $\esssup$ is the essential supremum operator. In the next section, we propose a stopping time for \eqref{eqn:lorden_formulation}.
	
	A closely related topic is \gls{TCD}\cite{bojdecki1980probability,pollak2013shewhart,moustakides2014multiple,poor18FMA} where the change only occurs for a finite period of time and the objective is to detect if such a change has occurred within a predefined window or not instead of detecting the change as quickly as possible. 
	There are two widely adopted methods for the TCD problem, the window-limited CuSum stopping time \cite{guepie2012sequential} and the  \gls{FMA} stopping time\cite{poor18FMA}. The \gls{FMA} stopping time has been shown to perform well for the TCD problem, and we will use the FMA stopping time as a comparison in \cref{sec:numerical}.
	When $\nu_c < \nu_n$, our system model can be seen to be a generalization of the TCD problem variant where one seeks to detect the transient change as quickly as possible by letting $g_n = f$. In \cref{eqn:WSGLR-test-statistic} below, we propose a test statistic and stopping time for \cref{eqn:lorden_formulation}. By setting $g_n=f$ and $f_n=f$ in our test statistic, our proposed stopping time reduces to the window-limited CuSum stopping time with pre-change distribution $f$ and post-change distribution $g$. 
	
	\section{Test Statistic for QCD with Nuisance Change}\label{sec:teststats}

	In this section, we derive a test-statistic and stopping time for QCD under a nuisance change. Suppose that we observe the sequence $X_1, X_2,\ldots$ and know \emph{a priori} that the nuisance change does not take place (i.e., $\nu_n=\infty$), then Page's CuSum test statistic\cite{page54} given as
	$
	S_{\text{CuSum}}(t)=\max_{1\leq k\leq t+1}\sum_{i=k}^t \log\tfrac{g(X_i)}{f(X_i)},
	$
	can be used and we declare that a critical change has taken place at
	\begin{align}\label{eqn:rpt_sprt}
	\tau_{\text{CuSum}}(b) & =\inf\left\{t\ :\ S_{\text{CuSum}}(t)\geq b\right\}                                            \\
	& =\inf\left\{t\ :\ \max_{1\leq k\leq t+1}\sum_{i=k}^t \log \tfrac{g(X_i)}{f(X_i)}\geq b\right\},
	\end{align}
	where $b$ is a pre-determined threshold.
	The CuSum test statistics has a convenient recursion
	$
	S_{\text{CuSum}}(t+1)=\max\left\{S_{\text{CuSum}}(t)+\log\tfrac{g(X_t)}{f(X_t)} ,0\right\},
	$
	which allows the CuSum stopping time to be implemented efficiently.
	
	If the nuisance change takes places at a time $\nu_n<\infty$ and $\nu_n$ is known, a modification of Page's test statistic gives the following:
	\begin{align}
	S_{\text{CuSum}}(t)    & =\max_{1\leq k\leq t+1}\sum_{i=k}^t \log\tfrac{h_{\nu_n,1,i}(X_i)}{h_{\nu_n,\infty,i}(X_i)},                                \\
	\tau_{\text{CuSum}}(b) & =\inf\left\{t\ :\ \max_{1\leq k\leq t+1}\sum_{i=k}^t \log\tfrac{h_{\nu_n,1,i}(X_i)}{h_{\nu_n,\infty,i}(X_i)}\geq b\right\},
	\end{align}
	where $h_{\nu_n,1,i}(x)$ and $h_{\nu_n,\infty,i}(x)$ are as defined in \eqref{eqn:definition_of_h} and are the probability distributions corresponding to the cases where the critical change has already occurred or will never occur, respectively. Similar to the case where $\nu_n=\infty$, the CuSum test statistics admits a convenient recursion for efficient implementation. Furthermore, for both the cases mentioned above, $\tau_{\text{CuSum}}$ was shown to be asymptotically optimal by \cite{lorden71}.
	
	A naive approach is to utilize four variants of $\tau_{\text{CuSum}}$, one for detecting for a change in each of the cases: from $f$ to $f_n$, from $f$ to $g$, from $f$ to $g_n$, and from $f_n$ to $g_n$. In the first stage, we monitor for changes from $f$ to either $f_n$, $g$ or $g_n$. If a change to $f_n$ is detected, then we monitor for a change from $f_n$ to $g_n$. The difficulty in such an approach is that any false alarm or miss detection in the first stage propagates to the second stage. We demonstrate that such an approach is suboptimal in \cref{subsec:numerical_results_sim_signals}.
	
	In our problem formulation, the nuisance change-point $\nu_n$ is unknown. Replacing $\nu_n$ with its maximum likelihood estimator in both the numerator and denominator, we obtain the following \gls{GLR} test statistic and stopping time:
	\begin{align}
	\Lambda_{\text{GLR}}(k,t) & =\frac{\max_{k\leq j\leq t+1}\prod_{i=k}^th_{j,1,i}(X_i)}{\max_{k\leq j\leq t+1}\prod_{i=k}^th_{j,\infty,i}(X_i)}\nn
	& =\frac{\max_{k\leq j\leq t+1}\prod_{i=k}^{j-1}g(X_i)\prod_{i=j}^{t}g_n(X_i)}{\max_{k\leq j\leq t+1}\prod_{i=k}^{j-1}f(X_i)\prod_{i=j}^{t}f_n(X_i)}, \label{Lambda_GLR} \\
	S_{\text{GLR}}(t)         
	& =\max_{1\leq k\leq t+1} \log\Lambda_{\text{GLR}}(k,t),\\
	\tau_{\text{GLR}}(b) & =\inf\left\{t\ :\ S_{\text{GLR}}(t)\geq b\right\}. \label{tau_GLR}
	\end{align}
	From our simulations in \cref{subsec:numerical_results_sim_signals}, it turns out that \cref{tau_GLR} does not achieve the best trade-off between \gls{ADD} and \gls{ARL} to false alarm over a wide range of threshold values $b$. Furthermore, its ARL is challenging to characterize theoretically since the GLR test statistic $\Lambda_{\text{GLR}}(k,t)$ is not a likelihood ratio and standard techniques in the QCD literature (e.g., Theorem 6.16 of \cite{poor2009quickest}) cannot be used to analyze its ARL. This is a critical problem for practical applications that require us to pre-determine a suitable threshold $b$ to achieve a desired ARL.
	
	To develop a stopping time with ARL that can be characterized theoretically, we simplify the maximum likelihood estimation in the numerator of \cref{Lambda_GLR} to be the maximum of only two cases $j=k$ and $j=t+1$. This gives us the Simplified GLR (SGLR) test statistic and stopping time as follows:
	\begin{align}
	\Lambda_{\text{SGLR}}(k,t) & =\frac{\max\left\{\prod_{i=k}^{t}g(X_i),\prod_{i=k}^{t}g_n(X_i)\right\}}{\max_{k\leq j\leq t+1}\prod_{i=k}^th_{j,\infty,i}(X_i)} \label{eqn:SGLR_teststats} \\
	S_{\text{SGLR}}(t)         & =\max_{1\leq k\leq t+1} \log\Lambda_{\text{SGLR}}(k,t),                                                                          \\
	\tau_{\text{SGLR}}(b)      & =\inf\left\{t\ :\ S_{\text{SGLR}}(t)\geq b\right\}. \label{tau_SGLR}
	\end{align}
	
	Unlike the CuSum test statistic, the SGLR test statistic does not have a convenient recursion. Any implementation of the SGLR stopping time would require computational resources that increases with the number of samples observed.  The requirement on computational resources would be a significant limitation for many practical applications. To limit the computational resources required,in the same spirit as \cite{lai98}, we propose the Window-Limited SGLR (W-SGLR) test statistic and stopping time as follows:
	\begin{align}
	S_{\text{W-SGLR}}(t) & =\max_{t-m_b\leq k\leq t+1} \log\Lambda_{\text{SGLR}}(k,t),       \label{eqn:WSGLR-test-statistic}     \\
	\tauWSGLR(b)         & =\inf\left\{t\ :\ S_{\text{W-SGLR}}(t)\geq b\right\},\label{tau_WSGLR}
	\end{align}
	where the window size $m_b$ is chosen such that
	\begin{align}\label{mb}
	\liminf_{b\to\infty} \frac{m_b}{b}>I^{-1} \text{ and } \log m_b=o(b),
	\end{align}
	with
	\begin{align}
	\begin{split}
	I=\min\left\{\KLD{g}{f},\ \KLD{g}{f_n},\ \KLD{g_n}{f},\ \KLD{g_n}{f_n}\right\}, \label{I}
	\end{split}
	\end{align}
	and $o(b)$ denoting a term that goes to zero as $b\to\infty$. Window-limited test statistics were first introduced by \cite{willsky1976generalized}. The paper \cite{lai98} further discussed their properties and the choice of window size and thresholds. We make the following assumption.
	
	\begin{Assumption}\label{assumpt:moments}
		The first four moments of $\log\tfrac{f_n(X)}{f(X)}$ \gls{wrt} both $g$ and $g_n$ are finite, and $\rho_g, \rho_{g_n}\neq 0$, where we define
		$ \rho_{g}=\E{g}[\log \tfrac{f_n(X)}{f(X)}]$, $\sigma_{g}^2=\E{g}[\left(\log \tfrac{f_n(X)}{f(X)}-\rho_g\right)^2]$, $\omega_{g}^4=\E{g}[\left(\log \tfrac{f_n(X)}{f(X)}-\rho_g\right)^4]$, $\rho_{g_n}=\E{g_n}[\log \tfrac{f_n(X)}{f(X)}]$, $\sigma_{g_n}^2=\E{g_n}[\left(\log \tfrac{f_n(X)}{f(X)}-\rho_{g_n}\right)^2]$ and $\omega_{g_n}^4=\E{g_n}[\left(\log \tfrac{f_n(X)}{f(X)}-\rho_{g_n}\right)^4]$.
	\end{Assumption}
	
	In \cref{thm:main_result} of \cref{subsec:properties_stoptime}, we show that the proposed $\tauWSGLR(b)$ is asymptotically optimal as $b\to\infty$ under \cref{assumpt:moments} and an additional technical assumption. To do that, we first analyze the asymptotic properties of $\tauWSGLR$. We let
	\begin{align}
	\Lambda(k,t)   & =\tfrac{\prod_{i=k}^t g(X_i)}{\max_{k\leq j\leq t+1}\prod_{i=k}^{j-1}f(X_i)\prod_{i=j}^{t}f_n(X_i)},\label{Lambda}    \\
	\Lambda_n(k,t) & =\tfrac{\prod_{i=k}^t g_n(X_i)}{\max_{k\leq j\leq t+1}\prod_{i=k}^{j-1}f(X_i)\prod_{i=j}^{t}f_n(X_i)},\label{Lambdan}
	\end{align}
	and study their properties in \cref{subsec:properties_stats}. Then, using the relationships
	\begin{align}
	\Lambda_{\text{SGLR}}(k,t) & =\max\left\{\Lambda(k,t),\Lambda_n(k,t)\right\},                       \\
	\tau_{\text{SGLR}}(b)      & =\min\{\tau(b),\tau_n(b)\}, \label{tau_SGLR2}                          \\
	\tau_{\text{W-SGLR}}(b)    & =\min\{\widetilde{\tau}(b),\widetilde{\tau_n}(b)\}, \label{tau_WSGLR2}
	\end{align}
	where
	\begin{align}
	\begin{split}\tau(b)             & =\inf\{t\ :\ \max_{k\leq t}\log\Lambda(k,t)\geq b \},\\
	\tau_n(b)&=\inf\{t\ :\ \max_{k\leq t}\log\Lambda_n(k,t)\geq b \},                  \\
	\widetilde{\tau}(b) & =\inf\{t\ :\ \max_{t-m_b \leq k\leq t}\log\Lambda(k,t)\geq b \},\\
	\quad\widetilde{\tau}_n(b)&=\inf\{t\ :\ \max_{t-m_b \leq k\leq t}\log\Lambda_n(k,t)\geq b \},\end{split}\label{tau_taun} 
	\end{align}
	we finally show the asymptotic optimality of $\tauWSGLR$ under mild technical conditions in \cref{subsec:properties_stoptime}.

	\subsection{Log Likelihood Ratio Growth Rates}\label{subsec:properties_stats}
	
	In this subsection, we derive properties of $\Lambda$ and $\Lambda_n$ as defined in \cref{Lambda,Lambdan}, respectively. The stopping times $\tau_{\text{SGLR}}(b)$ and $\tau_{\text{W-SGLR}}(b)$ are defined by the first time the test statistics $S_{\text{SGLR}}$ and $S_{\text{W-SGLR}}$ cross the threshold $b$ respectively. The rates of growth, $\tfrac{1}{t-k+1}\log\Lambda(k,t)$ and $\tfrac{1}{t-k+1}\log\Lambda_n(k,t)$, allow us to understand the detection delay of these stopping times. We show that these rates of growth converge in probability as $t\to\infty$. In particular, the limit that the rate of growth converges to depends on the sign of $\rho_{g}$ and $\rho_{g_n}$.
	
	As the nuisance change point is unknown, the denominator of both $\Lambda$ and $\Lambda_n$ contains a maximization of the likelihood
	$
	\max_{k\leq j\leq t+1}\prod_{i=k}^{j-1}f(X_i)\prod_{i=j}^{t}f_n(X_i).
	$
	If the first moment $\rho_g<0$, the distribution $g$ is closer to the distribution $f$ as compared to $f_n$ in the \gls{KL} divergence sense. When the critical change point is at $\nu_c=1$ and no nuisance change has taken place, we expect the denominator to approach $\prod_{i=k}^{t}f(X_i)$. Thus, our statistic $\Lambda(k,t)$ can be approximated by $\prod_{i=k}^t\tfrac{g(X_i)}{f(X_i)}$. A similar argument can be made for $\Lambda_n$ when $\nu_n=\nu_c=1$. This observation is made precise in the following two propositions.
	
	\begin{Proposition}\label{prop:converge_to_post_nui_LR}
		Suppose that \cref{assumpt:moments} holds, and $\rho_{g}<0$. For any $\nu_c\leq k<\infty$ and $\epsilon>0$, we have
		\begin{align}
		& \lim_{t\to\infty}\P{\infty,\nu_c}(\left|\tfrac{\log \Lambda(k,t)}{t-k+1}-\tfrac{1}{t-k+1}\sum_{i=k}^t\log \tfrac{g(X_i)}{f(X_i)}\right|\geq\epsilon)=0,    \\
		& \lim_{t\to\infty}\P{\infty,\nu_c}(\left|\tfrac{\log\Lambda_n(k,t)}{t-k+1}-\tfrac{1}{t-k+1}\sum_{i=k}^t\log \tfrac{g_n(X_i)}{f(X_i)}\right|\geq\epsilon)=0.
		\end{align}
	\end{Proposition}
	\begin{IEEEproof}
		See \cref{sec:AppProp1}.
	\end{IEEEproof}

	\begin{Proposition}\label{prop:converge_to_pre_nui_LR}
		Suppose that \cref{assumpt:moments} holds, and $\rho_{g_n}>0$. For any $\nu_c\leq k<\infty$, $\nu_n<\infty$, and $\epsilon>0$, we have
		\begin{align}
		& \lim_{t\to\infty}\P{\nu_n,\nu_c}(\left|\tfrac{\log\Lambda(k,t)}{t-k+1}-\tfrac{1}{t-k+1}\sum_{i=k}^t\log \tfrac{g(X_i)}{f_n(X_i)}\right|\geq\epsilon)=0,     \\
		& \lim_{t\to\infty}\P{\nu_n,\nu_c}(\left|\tfrac{\log\Lambda_n(k,t)}{t-k+1}-\tfrac{1}{t-k+1}\sum_{i=k}^t\log \tfrac{g_n(X_i)}{f_n(X_i)}\right|\geq\epsilon)=0.
		\end{align}
	\end{Proposition}
	\begin{IEEEproof}
		See \cref{sec:AppProp1}.
	\end{IEEEproof}
	
	Using \cref{prop:converge_to_post_nui_LR,prop:converge_to_pre_nui_LR} together with the weak law of large numbers, we obtain the following result.
	\begin{Theorem}\label{thm:convergence_in_probability_special}
		Suppose that \cref{assumpt:moments} holds, $\rho_{g}=\KLD{g}{f}-\KLD{g}{f_n}<0$, and $\rho_{g_n}=\KLD{g_n}{f}-\KLD{g_n}{f_n}>0$. For any $\nu_c\leq k<\infty$,
		\begin{align*}
		\tfrac{\log \Lambda(k,t)}{t-k+1}   & \xrightarrow{\P{\infty,\nu_c}} \KLD{g}{f},   \\
		\tfrac{\log \Lambda_n(k,t)}{t-k+1} & \xrightarrow{\P{\infty,\nu_c}} \KLD{g}{f}-\KLD{g}{g_n},
		\end{align*}
		as $t\to\infty$. Furthermore, for any $\nu_n < \infty$,
		\begin{align*}
		\tfrac{\log \Lambda(k,t)}{t-k+1}   & \xrightarrow{\P{\nu_n,\nu_c}} \KLD{g_n}{f_n}-\KLD{g_n}{g},   \\
		\tfrac{\log \Lambda_n(k,t)}{t-k+1} & \xrightarrow{\P{\nu_n,\nu_c}} \KLD{g_n}{f_n},
		\end{align*}
		as $t\to\infty$.
	\end{Theorem}
	In \cref{thm:convergence_in_probability_special}, we have assumed that $\rho_{g}<0$ and $\rho_{g_n}>0$. If we vary the signs of $\rho_{g}$ and $\rho_{g_n}$, a similar argument to that provided in \cref{thm:convergence_in_probability_special} gives us the following result.
	\begin{Theorem}\label{thm:convergence_in_probability_general}
		Suppose that \cref{assumpt:moments} holds. For any $\nu_c\leq k<\infty$ and $\nu_n<\infty$, we have the following convergences in probability as $t\to\infty$ shown in the table below.
		
		\begin{center}
			\begin{tabular}{cc|c|c}
				\hline\hline
				\multicolumn{2}{c}{Sign of moment} & \multicolumn{1}{c}{under $\P{\infty,\nu_c}$} &\multicolumn{1}{c}{under $\P{\nu_n,\nu_c}$}                                                                                                                   \\
				\cline{1-2}\cline{3-3}\cline{4-4}
				$\rho_g$                           & $\rho_{g_n}$                                 & $\tfrac{\log \Lambda(k,t)}{t-k+1}$           & $\tfrac{\log \Lambda_n(k,t)}{t-k+1}$  \\
				\hline\hline
				>0                                 & >0                                           & $\KLD{g}{f_n}$         & $\KLD{g_n}{f_n}$ \\
				>0&<0& $\KLD{g}{f_n}$& $\KLD{g_n}{f}$\\
				<0                                 & >0                                           & $\KLD{g}{f}$   & $\KLD{g_n}{f_n}$ \\
				<0                                 & <0                                           & $\KLD{g}{f}$    & $\KLD{g_n}{f}$   \\
				\hline
			\end{tabular}
		\end{center}
	\end{Theorem}
	\Cref{thm:convergence_in_probability_general} gives us the average rate of growth of the statistics $\log \Lambda(k,t)$ and $\log\Lambda_n(k,t)$. Since $I$ in \cref{I} is the minimum of the growth rates in \cref{thm:convergence_in_probability_general}, we see that the average growth rate of these statistics is at least $I$ regardless of the signs of $\rho_g$ and $\rho_{g_n}$. This suggests that the WADD of $\tau_{\text{W-SGLR}}(b)$ grows linearly with respect to $b$ with a gradient bounded above by $I$. This observation is made precise in the next subsection.
	
	\subsection{Conditions for Asymptotic Optimality}\label{subsec:properties_stoptime}
	In this subsection, we establish the asymptotic \gls{WADD}-\gls{ARL} trade-off under \cref{assumpt:moments} and provide a sufficient condition for $\tau_{\text{W-SGLR}}$ to be asymptotically optimal. In particular, we show that $\tau_{\text{W-SGLR}}$ is asymptotically optimal if in addition to \cref{assumpt:moments}, the following assumption holds.
	
	\begin{Assumption}\label{assumpt:kldiv}
		The KL divergences $\KLD{g}{f}$, $\KLD{g}{f_n}$, $\KLD{g_n}{f}$ and $\KLD{g_n}{f_n}$ satisfy
		\begin{align*}
		\KLD{g}{f_n}>\min\left\{\KLD{g}{f},\ \KLD{g_n}{f},\ \KLD{g_n}{f_n}\right\}.  
		\end{align*}		
	\end{Assumption}
	\cref{assumpt:kldiv} essentially says that $g$ cannot be too similar to $f_n$, which makes intuitive sense as otherwise it is difficult to distinguish the critical change from the nuisance change (see \cref{fig:h_model}). A sufficient condition for \cref{assumpt:kldiv} is $\rho_g <0$ as assumed in \cref{thm:convergence_in_probability_special}. For example, in the problem of spectrum sensing in cognitive radio\cite{poor08qcdspectrum}, we are often interested in detecting a variance change of a signal generated by independent Gaussian distributions. Furthermore, in many signal processing applications, a change in mean may be due to sensor drift as a result of long duration monitoring. This change in mean is usually not of interest and interferes with the actual signal processing task\cite{Wang17drift,li2015Driftdetection,Takruri07drift}. A typical signal model of this type is given by $f=\calN(\mu_0,\sigma_0^2)$, $f_n=\calN(\mu_1,\sigma_0^2)$, $g=\calN(\mu_0,\sigma_1^2)$, and $g_n=\calN(\mu_1,\sigma_1^2)$ with $\mu_0,\mu_1\in\mathbb{R}$, $\mu_0\neq\mu_1$, $\sigma_0,\sigma_1\in\mathbb{R}_{>0}$. While \cref{assumpt:kldiv} may seem artificial, it is shown in \cref{lem:sufficient_conditions_exponential_family} that this model satisfies $\rho_g <0$ and hence \cref{assumpt:kldiv}. In this case, the W-SGLR stopping time achieves asymptotic optimality.

	We use techniques from the proof of Theorem 6.16 in \cite{poor2009quickest} to obtain a lower bound for the ARL of $\tau_{\text{SGLR}}$ in \cref{tau_SGLR}. Since $\tauWSGLR\geq \tau_{\text{SGLR}}$, the same lower bound also applies for the ARL of $\tau_{\text{W-SGLR}}$ in \cref{tau_WSGLR}. In the previous subsection, we have shown that the rate of growth of the statistics $\Lambda$ and $\Lambda_n$ converge to constants as $t\to\infty$. This means that, asymptotically, $\Lambda$ and $\Lambda_n$ grow linearly \gls{wrt} $t$. Heuristically, this implies that the WADD of the stopping times, $\tau(b)$ and $\tau_n(b)$, grow linearly \gls{wrt} the threshold $b$ in \cref{tau_taun}.
	
	\Cref{lem:fa_prob} derives an upper-bound for the probability of a false alarm for a stopping time related to $\tauWSGLR(b)$. Following Theorem 6.16 in \cite{poor2009quickest}, this upper bound then yields a lower bound for the ARL of $\tau_{\text{W-SGLR}}(b)$ in \cref{thm:main_result}.
	
	\begin{Lemma}\label{lem:fa_prob}
		For $b>0$, let $\eta^k,\eta^k_n$ be stopping times defined by
		\begin{align*}
		\eta^k   & =\inf\{t\geq k :\ \log\Lambda(k,t)\geq b\},      \\
		\eta^k_n & =\inf\{t\geq k :\ \log \Lambda_n(k,t)\geq b\},
		\end{align*}
		so that $\tau(b)=\inf_{k \geq 1} \eta^k$ and $\tau_n(b)=\inf_{k \geq 1} \eta^k_n$ in \cref{tau_taun}. For any $\nu_n\in\mathbb{N}\cup\{\infty\}$, we have
		\begin{align}
		\P{\nu_n,\infty}(\eta^1<\infty)\leq e^{-b}\quad \text{and}\quad\P{\nu_n,\infty}(\eta_n^1<\infty)\leq e^{-b},
		\end{align}
		and
		\begin{align}\label{ARL_tauWSGLR}
		\E{\nu_n,\infty}[\tauWSGLR(b)] \geq \E{\nu_n,\infty}[\tauSGLR(b)] \geq \ofrac{2}e^b.
		\end{align}
		
	\end{Lemma}
	\begin{IEEEproof}
		See \cref{sec:AppLem1}.
	\end{IEEEproof}
	
	The next lemma checks that our proposed stopping time satisfies the assumption required in \cite{lai98} to relate the asymptotic upper-bound for the WADD to the threshold $b$ in \cref{prop:ADD}.
	
	\begin{Lemma}\label{lem:limsup_assumption}
		Suppose that \cref{assumpt:moments} holds.	For any $\delta>0$, we have
		\begin{enumerate}[(i)]
			\item\label{eqn:limsup:it1}$\displaystyle\lim_{t\to\infty}\quad\ \ \sup_{{\nu_n\in\mathbb{N}, 1\leq \nu_c \leq k}} \P{\nu_n,\nu_c}(\tfrac{1}{t}\log \Lambda(k,k+t-1)-I\leq -\delta)=0$,\ and
			\item\label{eqn:limsup:it2}$\displaystyle\lim_{t\to\infty}\sup_{1\leq \nu_c \leq k} \P{\infty,\nu_c}(\tfrac{1}{t}\log \Lambda(k,k+t-1)-I\leq -\delta)=0$.
		\end{enumerate}
	\end{Lemma}
	\begin{IEEEproof}
		See \cref{sec:AppLem2}.
	\end{IEEEproof}
	
	\begin{Proposition}\label{prop:ADD}
		Suppose that \cref{assumpt:moments} holds. There exists a $B$ such that for all $b \geq B$, we have
		\begin{enumerate}[(i)]
			\item\label{prop:ADD:it1} $\displaystyle\sup_{\nu_n,\nu_c\geq1}\esssup \E{\nu_n,\nu_c}[(\widetilde{\tau}_n(b)-\nu_c+1)^+]{X_1,\ldots,X_{\nu_c-1}}\\\leq (I^{-1}+o(1))b$, and
			\item\label{prop:ADD:it2} $\displaystyle\sup_{\nu_c\geq1}\esssup \E{\infty,\nu_c}[(\widetilde{\tau}(b)-\nu_c+1)^+]{X_1,\ldots,X_{\nu_c-1}}\leq (I^{-1}+o(1))b$.
		\end{enumerate}
	\end{Proposition}
	\begin{IEEEproof}
		See \cref{sec:AppProp3}.
	\end{IEEEproof}
	
	Finally, we show the asymptotic optimality of $\tauWSGLR$ in the following result.
	
	\begin{Theorem}\label{thm:main_result}
		Suppose that \cref{assumpt:moments} holds. For any $b>0$,
		\begin{align}
		\ARL(\tauWSGLR(b))        & \geq \tfrac{1}{2}e^b, \\
		\text{WADD}(\tauWSGLR(b)) & \leq(I^{-1}+o(1))b, \label{WSGLR_opt}
		\end{align}
		where $o(1)$ is a term going to zero as $b\to\infty$. Furthermore, if \cref{assumpt:kldiv} holds, then the stopping time $\tauWSGLR(b)$ is asymptotically optimal for the problem \eqref{eqn:lorden_formulation} as $b\to\infty$.
	\end{Theorem}
	\begin{IEEEproof}
		See \cref{sec:AppThm1}.
	\end{IEEEproof}
	
	In \cref{thm:main_result}, we have shown that $\tauWSGLR$ is asymptotically optimal under \cref{assumpt:moments} and \cref{assumpt:kldiv}. In the next lemma, we derive sufficient conditions for \cref{assumpt:kldiv} when $f,f_n,g,g_n$ belong to an exponential family.
	
	\begin{Lemma}\label{lem:sufficient_conditions_exponential_family}
		Suppose that $f$, $f_n$, $g$, $g_n\in\{\phi :\ \phi(x)=h(x)\exp\left(\sum_{i=1}^s B_i(\theta)T_i(x)-A(\theta)\right)\}$, an exponential family of distributions on $\mathbb{R}^N$ with parameters $\theta=\theta_f,\theta_{f_n},\theta_{g},\theta_{g_n},$ respectively. Here, $T_i\in\mathbb{R}^N\times\mathbb{R}$ and $A, B_i\in\mathbb{R}^M\times\mathbb{R}$ for $i=1,\ldots,s$. If any of the following inequalities hold:
		\begin{align}
		&A(\theta_{f_n})-A(\theta_{f})-\sum_{i=1}^s\left(B_i(\theta_{f_n})-B_i(\theta_{f})\right)\E{g}[T_i(X)]>0, \label{eqn:condition_for_exp_fam}\\
		\begin{split}
		&A(\theta_{f_n})-A(\theta_{g})-\sum_{i=1}^s\left(B_i(\theta_{f_n})-B_i(\theta_{g})\right)\E{g}[T_i(X)]\\
		&>A(\theta_{f})-A(\theta_{g_n})-\sum_{i=1}^s\left(B_i(\theta_{f})-B_i(\theta_{g_n})\right)\E{g_n}[T_i(X)],
		\end{split}\label{eqn:condition_for_exp_fam_1}\\
		\begin{split}
		&A(\theta_{f_n})-A(\theta_{g})-\sum_{i=1}^s\left(B_i(\theta_{f_n})-B_i(\theta_{g})\right)\E{g}[T_i(X)]\\
		&>A(\theta_{f_n})-A(\theta_{g_n})-\sum_{i=1}^s\left(B_i(\theta_{f_n})-B_i(\theta_{g_n})\right)\E{g_n}[T_i(X)],
		\end{split}\label{eqn:condition_for_exp_fam_2}
		\end{align}
		then \cref{assumpt:kldiv} holds.
		
		In particular, if $f=\calN(\mu_0,\sigma_0^2), f_n=\calN(\mu_1,\sigma_0^2),g=\calN(\mu_0,\sigma_1^2)$, and $g_n=\calN(\mu_1,\sigma_1^2)$ with $\mu_0,\mu_1\in\mathbb{R}$, $\mu_0\neq\mu_1$, $\sigma_0,\sigma_1\in\mathbb{R}_{>0}$, and $\sigma_0\neq \sigma_1$, \cref{assumpt:kldiv} holds.
	\end{Lemma}
	
	\begin{IEEEproof}
		To show that \cref{eqn:condition_for_exp_fam} implies \cref{assumpt:kldiv}, we rearrange the terms on the \gls{LHS} of \cref{eqn:condition_for_exp_fam} to obtain
		\begin{align}
		\E{g}[\log\tfrac{f_n(X)}{f(X)}]=\E{g}[\log\tfrac{ h(X)\exp\left(\sum_{i=1}^{s}B_i(\theta_{f_n})T_i(X)-A(\theta_{f_n})\right)}{h(X)\exp\left(\sum_{i=1}^{s}B_i(\theta_{f})T_i(X)-A(\theta_{f})\right)}]< 0.\label{eqn:exp_fam_proof_1}
		\end{align}
		This implies that
		\begin{align}
		\KLD{g}{f}< \KLD{g}{f_n},\label{eqn:exp_fam_proof_2}
		\end{align}
		and hence \cref{assumpt:kldiv} holds. A similar argument shows that \cref{eqn:condition_for_exp_fam_1,eqn:condition_for_exp_fam_2} imply $\KLD{g_n}{f}<\KLD{g}{f_n}$ and $\KLD{g_n}{f_n}<\KLD{g}{f_n}$, respectively.

		If $f=\calN(\mu_0,\sigma_0^2),\  f_n=\calN(\mu_1,\sigma_0^2),\ g=\calN(\mu_0,\sigma_1^2)$, and $g_n=\calN(\mu_1,\sigma_1^2)$ with $\mu_0,\mu_1\in\mathbb{R}$, $\mu_0\neq\mu_1$, $\sigma_0,\sigma_1\in\mathbb{R}_{>0}$, and $\sigma_0\neq \sigma_1$, we can define $\theta_f=[\mu_0,\sigma_0^2],\ \theta_{f_n}=[\mu_1,\sigma_0^2],\ \theta_g=[\mu_0,\sigma_1^2],\ \theta_{g_n}=[\mu_1,\sigma_1^2]$, with the functions $B_1(\mu,\sigma^2)=\mu/\sigma^2$, $B_2(\mu,\sigma^2)=\tfrac{-1}{2\sigma^2}$, $T_1[X]=X$, $T_2[X]=X^2$ and $A(\mu,\sigma^2)=\tfrac{\mu^2}{2\sigma^2}+\log \sigma$. The \gls{LHS} of \cref{eqn:condition_for_exp_fam} becomes $\tfrac{\mu_1^2-\mu_0^2}{\sigma_0^2}-\left(\tfrac{{\mu_1-\mu_0}}{\sigma_0}\right)\mu_0$. Simplifying, the \gls{LHS} of \cref{eqn:condition_for_exp_fam} becomes $\tfrac{(\mu_1-\mu_0)^2}{2\sigma^2_0}$. Thus, for any $\mu_1\neq\mu_0$ and $\sigma_0,\sigma_1\in\mathbb{R}_{>0}$, the inequality \cref{eqn:condition_for_exp_fam} holds. The proof is now complete.
	\end{IEEEproof}

	\section{Parametrized Families of Post-Change Distributions}\label{sec:discussion}
	In many applications, the post-change distribution $g$ and nuisance post-change distribution $g_n$ may contain unknown parameters. In this section, we modify $\tau$ and $\tau_n$ in \cref{tau_taun} to obtain a \gls{GLRT}-based stopping time $\widehat{\tau}_{\text{W-SGLR}}$ for the following signal model: Let $\Theta\subseteq \Real^d$ be a set with non-empty interior and $X_1,X_2,\ldots$ be a sequence of independent random variables satisfying:
	$
	X_t \sim h_{\nu_n,\nu_c,\theta,\theta_n,t}
	$
	where \begin{align}
	h_{\nu_n,\nu_c,\theta,\theta_n,t}(\cdot)=\begin{cases}
	f(\cdot)            & \text{if $t<\min\{\nu_c,\nu_n\}$,}     \\
	f_n(\cdot)          & \text{if $\nu_n\leq t<\nu_c$,}         \\
	g(\cdot;\theta)     & \text{if $\nu_c\leq t<\nu_n$,}         \\
	g_n(\cdot;\theta_n) & \text{if $\max\{\nu_c,\nu_n\}\leq t$,} \\
	\end{cases}
	\end{align} and $\theta,\theta_n\in\text{Int}(\Theta)$, the interior of $\Theta$. We derive a lower bound for the ARL of $\widehat{\tau}_{\text{W-SGLR}}$ under the following assumption.
	
	\begin{Assumption}\label{assup:compact2nddiff}
		$\Theta$ is a compact $d$-dimensional sub-manifold of $\mathbb{R}^d$. The pdfs of the post-change distributions $g(\cdot;\theta)$ and nuisance post-change distribution $g_n(\cdot;\theta)$ are twice continuously differentiable \gls{wrt} $\theta$.
	\end{Assumption}
	
	A commonly used method to handle unknown parameters is to replace the likelihood ratio $\Lambda(k,t)$ with the generalized likelihood ratio. We define the generalized W-SGLR test statistic $\widehat{S}_{\text{W-SGLR}}$ as
	\begin{align}
	\widehat{\Lambda}_{\text{W-SGLR}}(k,t,\theta) & =\tfrac{\max\left\{\prod_{i=k}^{t}g(X_i;\theta),\prod_{i=k}^{t}g_n(X_i;\theta)\right\}}{\max_{k\leq j\leq t+1}\prod_{i=k}^th_{j,\infty,\theta,\theta,i}(X_i)},\nonumber \\
	\widehat{S}_{\text{W-SGLR}}(t) & =\max_{t-m_{b}+1\leq k\leq t-m_{b}'}\max_{\theta\in\Theta}\log\widehat{\Lambda}(k,t,\theta),\label{eqn:generalised_wsglr}
	\end{align} 
	where the minimal delay $m_b'$ is required to prevent difficulties of under-determination when performing maximum likelihood estimation of the parameter $\theta$. While \cref{eqn:generalised_wsglr} is commonly used, the maximization over $\Theta$ make it difficult to theoretically quantify the \gls{ARL} of the stopping time $\inf\{t:\ \widehat{S}_{\text{W-SGLR}}(t) \geq b\}$. To work around this problem, we modify the stopping times $\tau$ and $\tau_n$ as follows. Let $\lambda_{\max}(A)$ denote the largest eigenvalue of the symmetric matrix $A$. Fix $m_b'\geq 0$. We let
	\begin{align}
	& \widehat{\Lambda}(k,t,\theta)=\tfrac{\prod_{i=k}^{t}g(X_i;\theta)}{\max_{k\leq j\leq t+1}\prod_{i=k}^th_{j,\infty,\theta,\theta,i}(X_i)},\\
	& \widehat{\tau}(b)=\min_{m_{b}'\leq l\leq m_{b}} \widehat{\eta}_l(b), \label{eqn:modified_SGLR1}\\
	\begin{split}
	&\widehat{\eta}_l(b)=\inf\bigg\{t\geq l:\ \log\widehat{\Lambda}(t-l+1,t,\widehat{\theta})\geq b,\\
	&\quad \sup_{\|\theta-\widehat{\theta}\|<1/\sqrt{b}}\lambda_{\max}\left(-\nabla^2\log\widehat{\Lambda}(t-l+1,t,\theta)\right)\leq b,\\
	&\quad\quad\text{and} \ \widehat{\theta}=\argmax_\theta \log\widehat{\Lambda}(t-l+1,t,\theta)\in \text{Int}(\Theta)  \bigg\},
	\end{split}\\
	& \widehat{\Lambda}_n(k,t,\theta)=\tfrac{\prod_{i=k}^{t}g_n(X_i;\theta)}{\max_{k\leq j\leq t+1}\prod_{i=k}^th_{j,\infty,\theta,\theta,i}(X_i)},\\
	& \widehat{\tau}_n(b)=\min_{m_{b}'\leq l\leq m_{b}} \widehat{\eta}_{n,l}(b)\label{eqn:modified_SGLR2},\\
	\begin{split}
	&\widehat{\eta}_{n,l}(b)=\inf\bigg\{t\geq l:\ \log\widehat{\Lambda}_n(t-l+1,t,\widehat{\theta}_n)\geq b,\\
	&\sup_{\|\theta-\widehat{\theta}_n\|<1/\sqrt{b}}\lambda_{\max}\left(-\nabla^2\log\widehat{\Lambda}_n(t-l+1,t,\theta)\right)\leq b,\\
	&\quad\text{and}\  \widehat{\theta}_n=\argmax_\theta \log\widehat{\Lambda}_n(t-l+1,t,\theta)\in \text{Int}(\Theta) \bigg\}.
	\end{split}
	\end{align}
	We define the generalized W-SGLR stopping time as
	\begin{align}
	\widehat{\tau}_{\text{W-SGLR}}(b)=\min\{\widehat{\tau}(b),\widehat{\tau}_n(b)\}.
	\end{align}
	Note that $\widehat{\tau}_{\text{W-SGLR}}$ is a modification of $\tauWSGLR$ with additional conditions required for stopping.
	
	The paper \cite{willsky1976generalized} first introduced window-limited generalized detection rules. We compute the false alarm probability of $\widehat{\eta}_l$ and $\widehat{\eta}_{n,l}$ . We then use this false alarm probability to obtain a lower bound for the ARL of $\widehat{\tau}$ and $\widehat{\tau}_n$ in Proposition~\ref{prop:generalised_arl}.

	\begin{Lemma}\label{lem:error_prob_glrt}
		Suppose that \cref{assup:compact2nddiff} holds. Given any $0<\delta< 1$, there exists $b_\delta>0$ such that $\P{\nu_n,\infty}(\widehat{\eta}_k<\infty)\leq \exp\left(-(1-\delta)b\right)$
		and $\P{\nu_n,\infty}(\widehat{\eta}_{n,l}<\infty)\leq \exp\left(-(1-\delta)b\right)$ for any $\nu_n\in\mathbb{N}\cup\{\infty\}$ and $b\geq b_\delta$.
	\end{Lemma}
	\begin{IEEEproof}
		As the proof is similar to Lemma 2 in \cite{lai98}, we omit it here and refer the reader to the extended version in \cite{lau2019quickestarxiv}.
	\end{IEEEproof}

	\begin{Proposition}\label{prop:generalised_arl}
		Suppose that \cref{assup:compact2nddiff} holds. For any $0< \delta< 1$, there exists $b_\delta>0$ such that for all $b\geq b_\delta$ and $\nu_n\in\mathbb{N}\cup\{\infty\}$, we have
		\begin{align}
		\ARL(\widehat{\tau}_{\text{W-SGLR}}(b))=\inf_{\nu_n}\E_{\nu_n,\infty}[\widehat{\tau}_{\text{W-SGLR}}(b)]\geq \tfrac{1}{2}e^{(1-\delta)b}.
		\end{align}
	\end{Proposition}
	\begin{IEEEproof}
		Fix $0< \delta<1$. By Lemma \ref{lem:error_prob_glrt}, there exists $b_\delta\geq 0$ such that
		$
		\P{\nu_n,\infty}(\min\{\widehat{\eta}_1,\widehat{\eta}_{n,1}\}<\infty)\leq 2\exp(-(1-\delta)b)
		$ for all $b\geq b_\delta$. 
		Applying results from Theorem 6.16 in~\cite{poor2009quickest}, we obtain
		\begin{align}
		\E_{\nu_n,\infty}[\widehat{\tau}_{\text{W-SGLR}}]\geq \tfrac{1}{2}e^{(1-\delta)b}
		\end{align}
		for all $b\geq b_\delta$. Taking infimum over $\nu_n$, we have
		$
		\ARL(\widehat{\tau}_{\text{W-SGLR}}(b))\geq \tfrac{1}{2}e^{(1-\delta)b},
		$
		and the proof is complete.
	\end{IEEEproof}

	\section{Numerical Results}\label{sec:numerical}
	In this section, we first illustrate the performance of the proposed W-SGLR stopping time under the assumption that the distributions $f,f_n,g$ and $g_n$ are known. Next, we illustrate the performance of the proposed generalised W-SGLR stopping time when $g$ and $g_n$ belongs to a parametrized family of distributions. Finally, we evaluate the performance of the proposed W-SGLR stopping time on real data from the Case Western Reserve University Bearing Dataset\cite{smith2015rolling}.
	
	\subsection{W-SGLR on Synthetic Data Satisfying \texorpdfstring{\cref{assumpt:kldiv}}{\ref{assumpt:kldiv}}}\label{subsec:numerical_results_sim_signals}
	
	\begin{figure}[!htb]
		\begin{subfigure}[b]{\columnwidth}
			\centering
			\includegraphics[width=0.85\columnwidth]{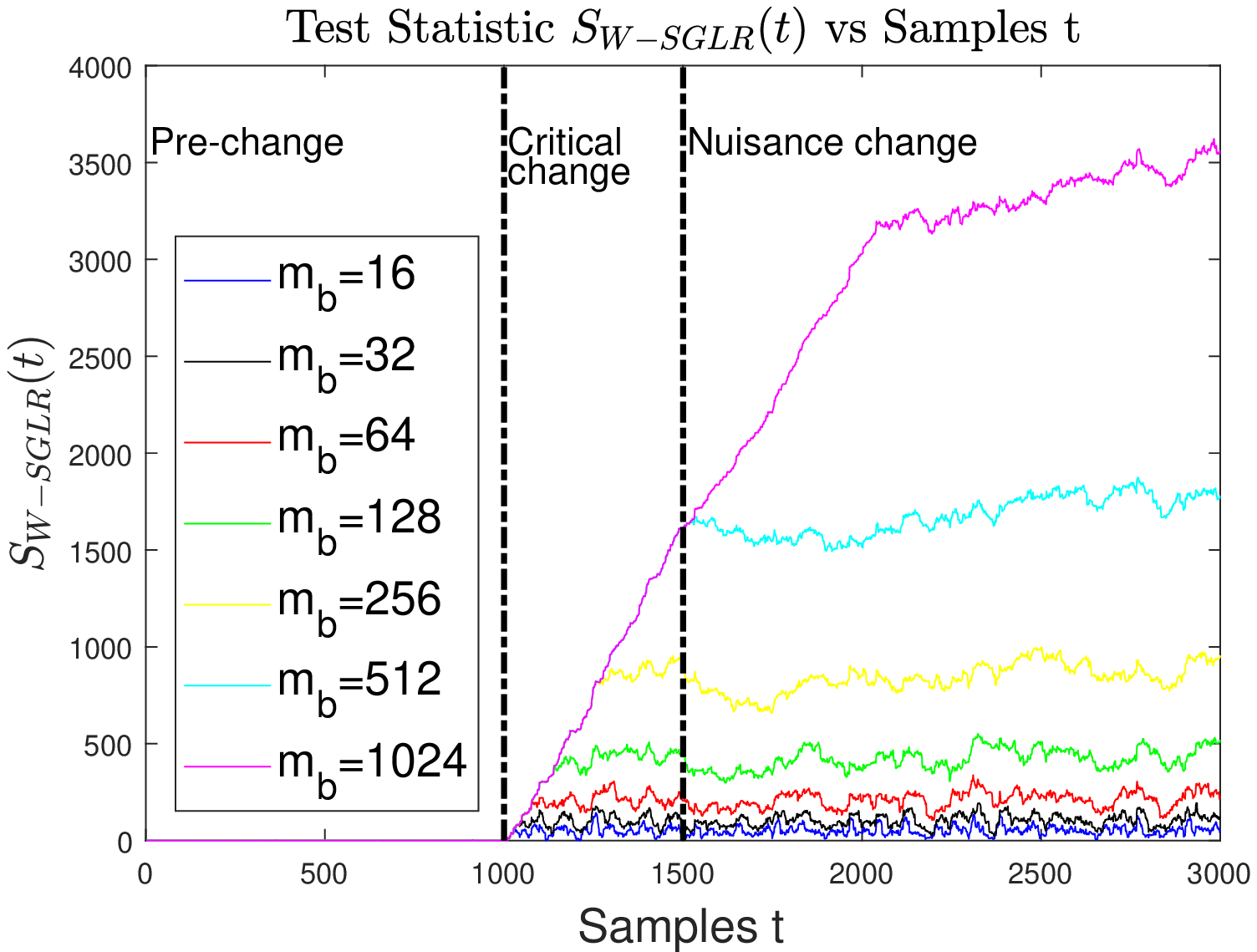}
			\caption{$\nu_c=1000$, $\nu_n=1500$}
			\label{fig:crit1000_nui1500}
		\end{subfigure}
		\begin{subfigure}[b]{\columnwidth}
			\centering
			\includegraphics[width=0.85\columnwidth]{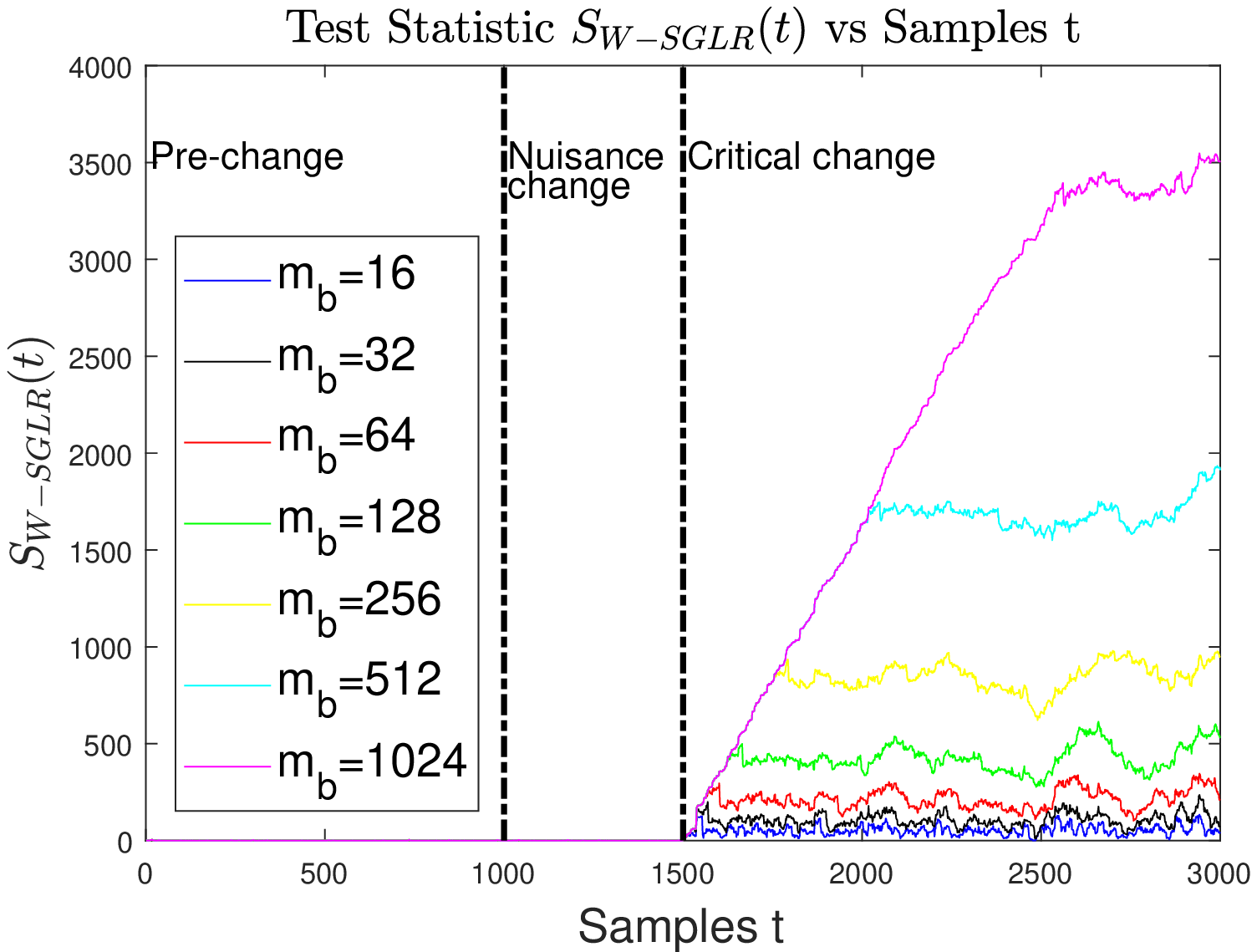}
			\caption{$\nu_c=1500$, $\nu_n=1000$}	\label{fig:crit1500_nui1000}
		\end{subfigure}
		\caption{W-SGLR test statistics $S_{\text{W-SGLR}}(t)$ with $f=\calN(0,1)$, $f_n=\calN(2,1)$, $g=\calN(0,10)$, $g_n=\calN(2,10)$ and $I=3.34$.}
		\label{fig:satassumpt}
	\end{figure}
	
	In our first set of simulations, we let $f=\calN(0,1)$,  $f_n=\calN(2,1)$, $g=\calN(0,10)$, and $g_n=\calN(2,10)$ where the critical change is a change in variance and the nuisance change is a change in mean (see example after \cref{assumpt:kldiv} for motivation).  We ran the simulations with two change-point configurations to illustrate the behaviour of the W-SGLR test statistic for different window-sizes. In \figref{fig:crit1000_nui1500}, we set $\nu_{c}=1000,\nu_n=1500$, while in \figref{fig:crit1500_nui1000}, we set $\nu_{c}=1500,\nu_n=1000$. In \figref{fig:crit1000_nui1500} and \figref{fig:crit1500_nui1000}, the test statistic $S_{\text{W-SGLR}}(t)$ remains low before the critical change-point and grows linearly with the gradient of at least $I=3.34$, as described in \cref{I}, after the critical change-point. This trend continues until the test statistic approximately achieves the value of $m_b I$. From our choice of $m_b$ in \cref{mb}, we see that $m_b I > b$ for $b$ large, i.e., our $\tauWSGLR$ is able to detect the critical change given sufficient delay for every choice of $b$ sufficiently large. However, choosing a larger $m_b$ is more resistant to outlier noise. For example, in \figref{fig:crit1000_nui1500}, when $m_b=1024$, we note that the test statistic continues to grow linearly with the gradient $I$ even after the nuisance change point. The trade-off is the increase in memory requirement and computational complexity. In \figref{fig:crit1500_nui1000}, we note that the test statistic continues to remain low during the period between the nuisance and the critical change point. This demonstrates that $\tauWSGLR$ is oblivious to the nuisance change.
	
	\begin{figure}[!htbp]
		\centering
		\includegraphics[width=0.85\columnwidth]{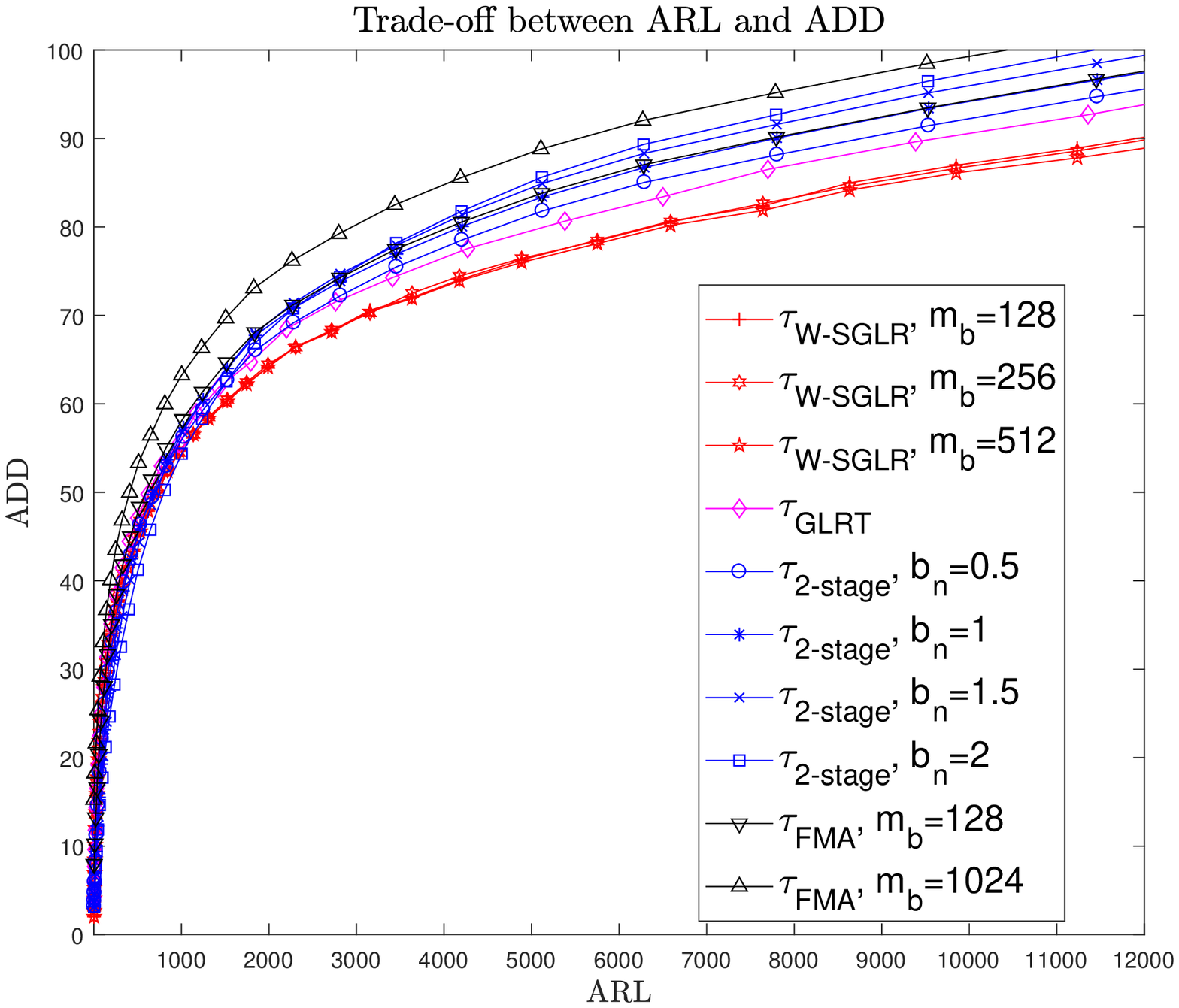}
		\caption{Comparison of the trade-off performance for $\tau_{\text{W-SGLR}}$, $\tau_{\text{GLRT}}$, $\tau_{\text{FMA}}$ and $\tau_{\text{2-stage}}$.}
		\label{fig:comparisonwithtwostage}
	\end{figure}
	
	Next, we compare $\tauWSGLR$, the GLRT stopping time $\tau_{\text{GLRT}}$ developed in \cite{lau2018quickest}, the finite-moving average (FMA) stopping time $\tau_{\text{FMA}}$ and a naive 2-stage CuSum stopping time denoted as $\tau_{\text{2-stage}}$. Following ideas from the \gls{TCD} literature\cite{poor18FMA}, the FMA stopping time $\tau_{\text{FMA}}$ is constructed by replacing the maximum in the test statistic \cref{eqn:WSGLR-test-statistic} with a sum across the entire window, i.e. setting $k=t-m_b$. 
	It should be noted that while the \gls{FMA} stopping time has been shown to perform well for the TCD problem, there are no guarantees that it will perform as well for the QCD problem.
	The naive stopping time $\tau_{\text{2-stage}}$ is constructed from stopping times based on the CuSum stopping time described in \eqref{eqn:rpt_sprt} with
	$
	\tau_{\{p\to q\}}(b)=\inf\left\{t\ :\ \max_{1\leq k\leq t+1}\prod_{i=k}^t \tfrac{p(x_i)}{q(x_i)}>e^b\right\}
	$
	for any pair of pdfs $p$ with $q\neq p$. We consider four stopping times: $\tau_{f\to f_n}(b_n)$, $\tau_{f\to g}(b_c)$, $\tau_{f\to g_n}(b_c)$, and $\tau_{f_n\to g_n}(b_c)$, where the threshold for declaring a critical change is $b_c$ and the threshold for declaring a nuisance change is $b_n$. In the first stage, we apply the stopping times $\tau_{f\to g}(b_c)$, $\tau_{f\to g_n}(b_c)$ and $\tau_{f\to f_n}(b_n)$ to the observations. If $\tau_{f\to g}(b_c)$ or $\tau_{f\to g_n}(b_c)$ stops the process before $\tau_{f\to f_n}(b_n)$, we declare that a critical change has occurred and set $\tau_{\text{2-stage}}=\min\{\tau_{f\to g}(b_c),\tau_{f\to g_n}(b_c)\}$. Otherwise, we apply $\tau_{f_n\to g_n}(b_c)$ to the rest of the observations after the stopping time $\tau_{f\to g_n}(b_n)$ and set $\tau_{\text{2-stage}}=\tau_{f_n\to g_n}(b_c)$. 
	
	In our simulations, our signal is generated using $f=\mathcal{N}(0,1)$, $g=\mathcal{N}(0.5,1)$, $f_n=\mathcal{N}(0,2)$, $g_n=\mathcal{N}(0.5,2)$. Here, the critical change is a change in mean from $0$ to $0.5$, and the nuisance change is a change in variance from $1$ to $2$. We generate a signal of length $2^{16}=65536$ and independently select the nuisance change point and critical change point with uniform probability on the $2^{16}$ possible data points. A total of $2^{12}=4096$ signals are generated. We compare the trade-off between the \gls{ADD} and the empirical \gls{ARL} of the proposed $\tauWSGLR$, $\tau_{\text{GLR}}$, $\tau_{\text{FMA}}$ and $\tau_{\text{2-stage}}$ in \figref{fig:comparisonwithtwostage}. We observe that our proposed $\tauWSGLR$ achieves a lower \gls{ADD} as compared to $\tau_{\text{FMA}}$, $\tau_{\text{GLR}}$ and $\tau_{\text{2-stage}}$ for large empirical \gls{ARL}.

	\begin{figure}[!htb]
		\begin{subfigure}[b]{\columnwidth}
			\centering
			\includegraphics[width=0.85\columnwidth]{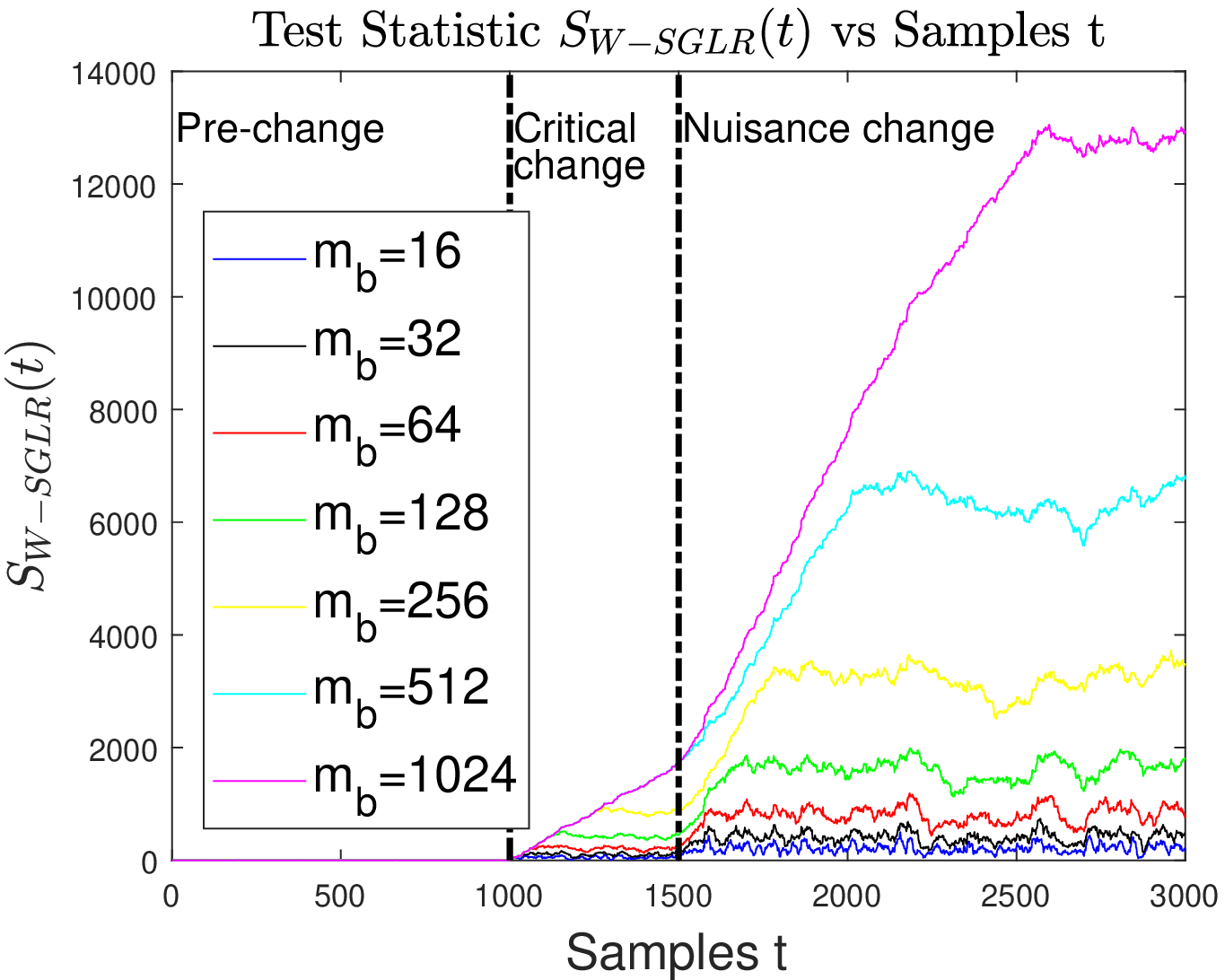}
			\caption{$\nu_c=1000$, $\nu_n=1500$}
			\label{fig:crit1000_nui1500_rho_negative}
		\end{subfigure}
		\begin{subfigure}[b]{\columnwidth}
			\centering
			\includegraphics[width=0.85\columnwidth]{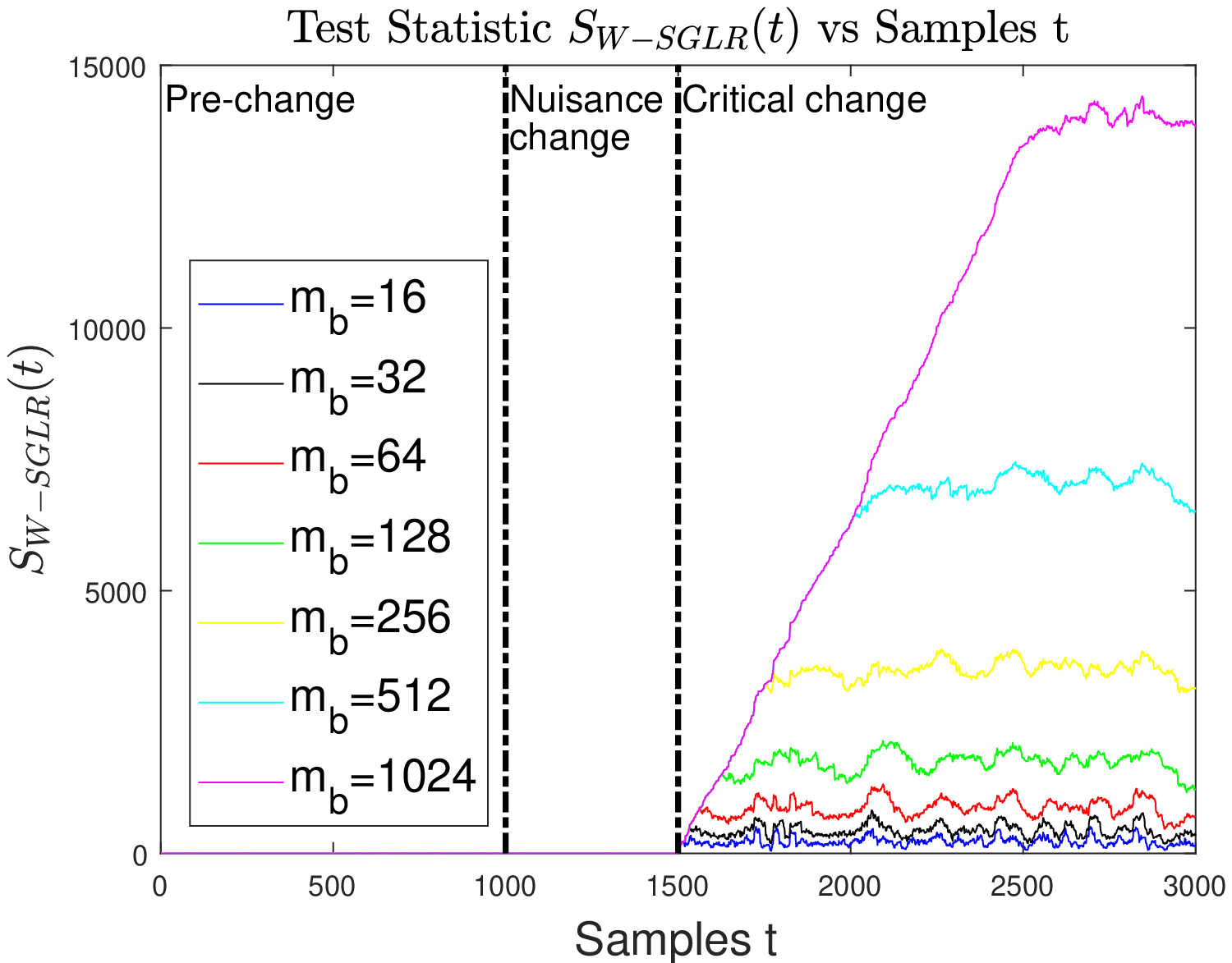}
			\caption{$\nu_c=1500$, $\nu_n=1000$}	\label{fig:crit1500_nui1000_rho_negative}
		\end{subfigure}
		\caption{W-SGLR test statistics $S_{\text{W-SGLR}}(t)$ with $f=\calN(0,1)$, $f_n=\calN(2,1)$, $g=\calN(0,10)$, $g_n=\calN(5,20)$ and $I=3.34$. }	\label{fig:rho_negative}
	\end{figure}

	In the next set of simulations, we let $f=\calN(0,1)$, $f_n=\calN(2,1)$, $g=\calN(0,10)$, and $g_n=\calN(5,20)$ where unlike the first set of simulations, the change in mean before and after the critical change point differs, and the change in variance before and after the nuisance change point differs. In \figref{fig:crit1000_nui1500_rho_negative}, we set $\nu_{c}=1000 < \nu_n=1500$. We see that the test statistic $S_{\text{W-SGLR}}(t)$ remains low before the critical change-point and grows linearly with a gradient of at least $I=3.34$ (cf.\ \cref{I}), after the critical change-point. This trend continues until the nuisance change point $\nu_n$ where rate of growth changes to $\E_g[\log\frac{g_n}{f_n}]>I$ until it approximately achieves the value of $m_b \E_g[\log\frac{g_n}{f_n}]$. While the growth of the test statistic after the change-point is not linear, the observation that the overall rate of growth from the critical change point is at least $I$ is consistent with \cref{lem:limsup_assumption}.  In \figref{fig:crit1500_nui1000_rho_negative}, we set $\nu_{c}=1500 > \nu_n=1000$ and note that the test statistic continues to remain low during the period between the nuisance and the critical change points. This demonstrates that $\tauWSGLR$ is oblivious to the nuisance change prior to the critical change-point.
	
	\subsection{W-SGLR on Synthetic Data Violating \texorpdfstring{\cref{assumpt:kldiv}}{\ref{assumpt:kldiv}}}\label{subsec:numerical_results_sim_signals_no_assumpt_2}
	
	\begin{figure}[!htb]
		\begin{subfigure}[b]{\columnwidth}
			\centering
			\includegraphics[width=0.85\columnwidth]{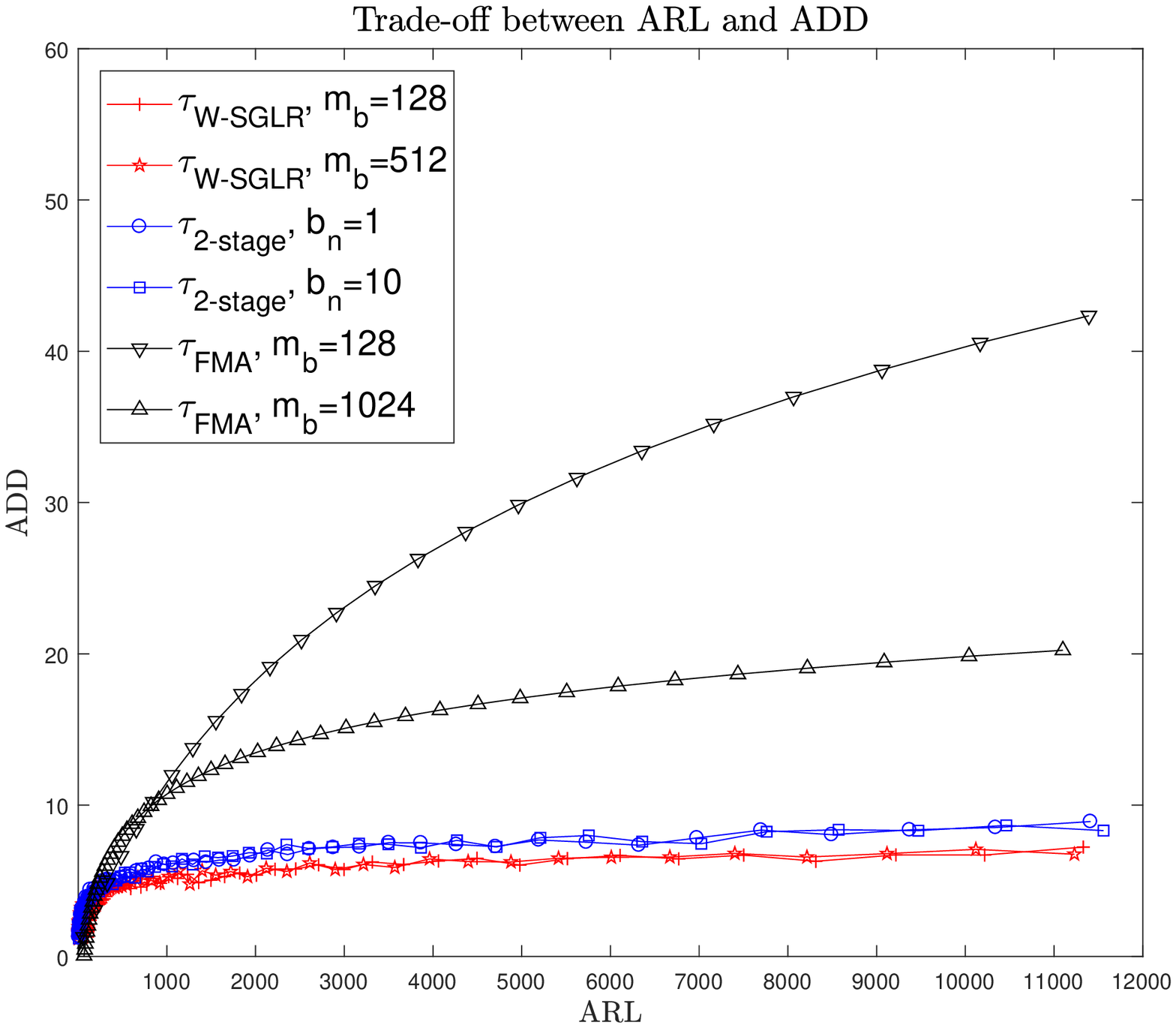}
			\caption{$\nu_n=2000$ and $\nu_c=4000$.}
			\label{fig:nui2000_crit4000_ARL_ADD}
		\end{subfigure}
		\ \\
		\begin{subfigure}[b]{\columnwidth}
			\centering
			\includegraphics[width=0.85\columnwidth]{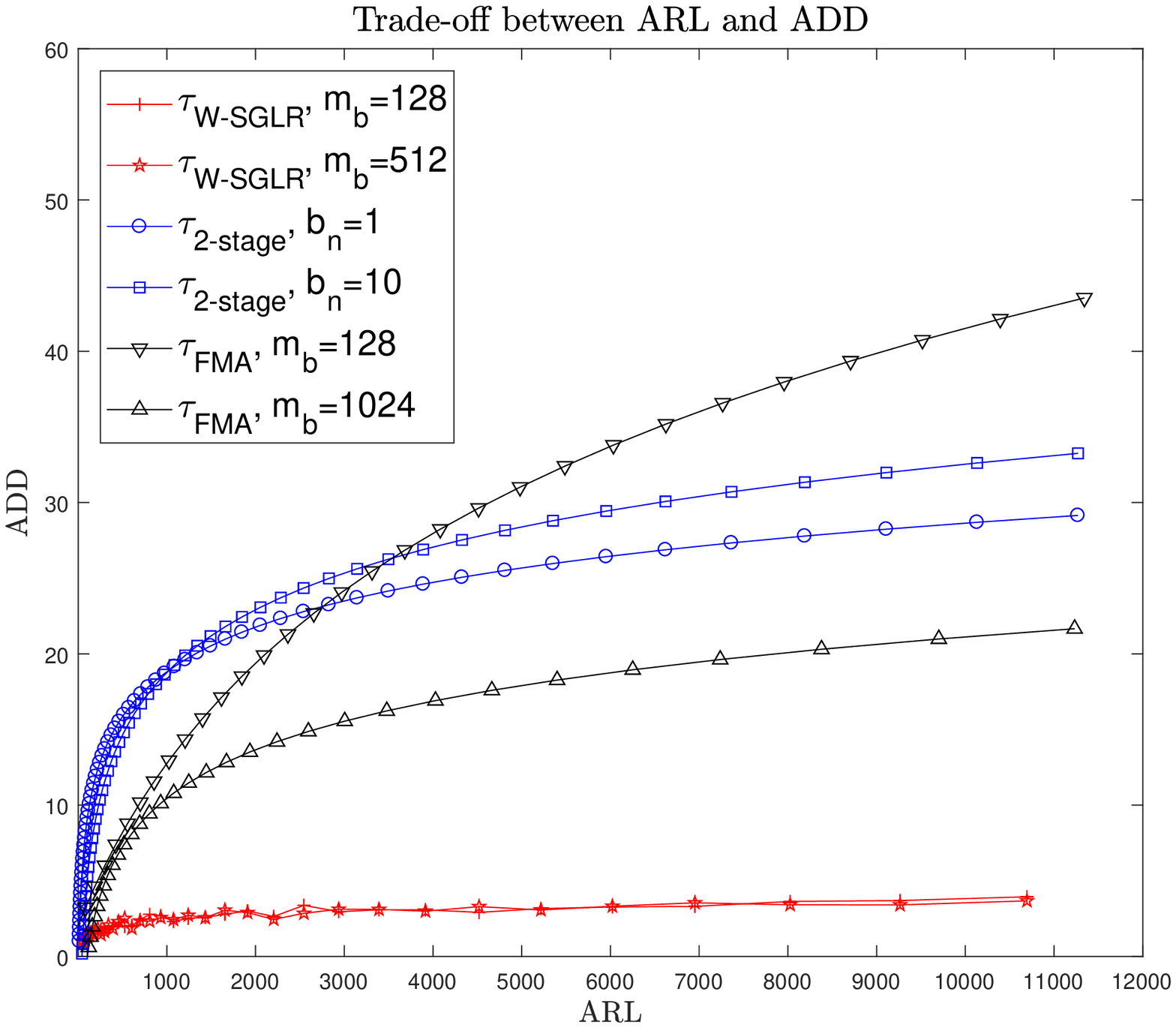}
			\caption{$\nu_n=4000$ and $\nu_c=2000$.}
			\label{fig:nui4000_crit2000_ARL_ADD}
		\end{subfigure}
		\caption{Comparison of the trade-off performance for $\tau_{\text{W-SGLR}}$, $\tau_{\text{2-stage}}$ and $\tau_{\text{FMA}}$ when \cref{assumpt:kldiv} is violated with $f=\calN(0,1)$, $f_n=\calN(2,5)$, $g=\calN(3,10)$ and $g_n=\calN(5,10)$.}\label{fig:non_asump_2_ARL_ADD}
	\end{figure}
	
	When \cref{assumpt:kldiv} is violated, \cref{thm:main_result} still provides the asymptotic trade-off between the ARL and the WADD. However, the asymptotic optimality of the W-SGLR stopping is not guaranteed. Here, we provide discussions and numerical simulations that suggests that the W-SGLR stopping time out-performs the two-stage stopping time and FMA stopping time with respect to \cref{eqn:lorden_formulation}.
	
	If \cref{assumpt:kldiv} is violated, from \cref{thm:main_result}, we have
	\begin{align*}
	\text{WADD}(\tauWSGLR(b)) & \leq\left(\KLD{g}{f_n}^{-1}+o(1)\right)b.
	\end{align*}
	This worst-case performance of the W-SGLR stopping time is achieved when $\nu_n=\infty$. For the rest of this discussion, we let $\nu_n=\infty$ to compare the two-stage stopping time with our proposed W-SGLR stopping time under this worst-case scenario. 
	
	Since $\KLD{g}{f}\geq\KLD{g}{f_n}$, the CuSum $$\max_{1\leq k\leq t+1}\sum_{i=k}^t\log\frac{f_n(X_i)}{f(X_i)}$$ associated with the stopping time $\tau_{f\to f_n}$ experiences a positive drift when $\nu_c\leq t<\nu_n$. Thus, for any finite threshold $b_n$ and sufficiently large $b_c$, the two-stage stopping time declares that a nuisance change has taken place and transits into the second stage after the critical change point. The CuSum associated with the stopping time $\tau_{f_n\to g_n}$ in the second stage is expected to grow at a rate of $\E{g}[\log \frac{g_n}{f_n}]$ for when $\nu_c\leq t<\nu_n$. In contrast, the W-SGLR test statistic, from \cref{thm:convergence_in_probability_general}, is expected to grow at a rate of $\E{g}[\log \frac{g}{f_n}]\geq \E{g}[\log \frac{g_n}{f_n}]$. Heuristically, this means that, when $\nu_n=\infty$, the $\WADD(\tau_{\text{2-stage}})\geq \WADD(\tau_{\text{W-SGLR}})$ as the $\ARL\to\infty$. It should also be noted that it is possible that the stopping time $\tau_{\text{2-stage}}$ fails completely when $\E{g}[\log \frac{g_n}{f_n}]$ is negative and $\nu_n=\infty$.
	
	In \cref{fig:nui4000_crit2000_ARL_ADD,fig:nui2000_crit4000_ARL_ADD}, we compare the trade-off between the ADD and the ARL of the different stopping times when \cref{assumpt:kldiv} is violated under the cases $(\nu_n,\nu_c)=(2000,4000)$ and $(\nu_n,\nu_c)=(4000,2000)$, respectively. To estimate the empirical ARL, the stopping times are applied to a set of $4096$ signals each of length $2^{16}$ with nuisance change point independently selected with uniform probability on the $2^{16}$ possible data points. To compute the corresponding ADD for the stopping times, they are applied to a set of $4096$ signals of length $4500$. It can be seen from both \cref{fig:nui4000_crit2000_ARL_ADD,fig:nui2000_crit4000_ARL_ADD} that the W-SGLR stopping time achieves a lower ADD as compared to both $\tau_{\text{2-stage}}$ and $\tau_{\text{FMA}}$ for large empirical \gls{ARL}. Consistent with our intuition, it can be seen that $\tau_{\text{W-SGLR}}$ significantly outperforms $\tau_{\text{2-stage}}$ when $\nu_c<\nu_n$.
	
	\subsection{Parametrized Post-Change Distributions}\label{subsec:parameterized}
	\begin{figure}[!htbp]
		\begin{subfigure}[b]{\columnwidth}
			\centering
			\includegraphics[width=0.85\columnwidth]{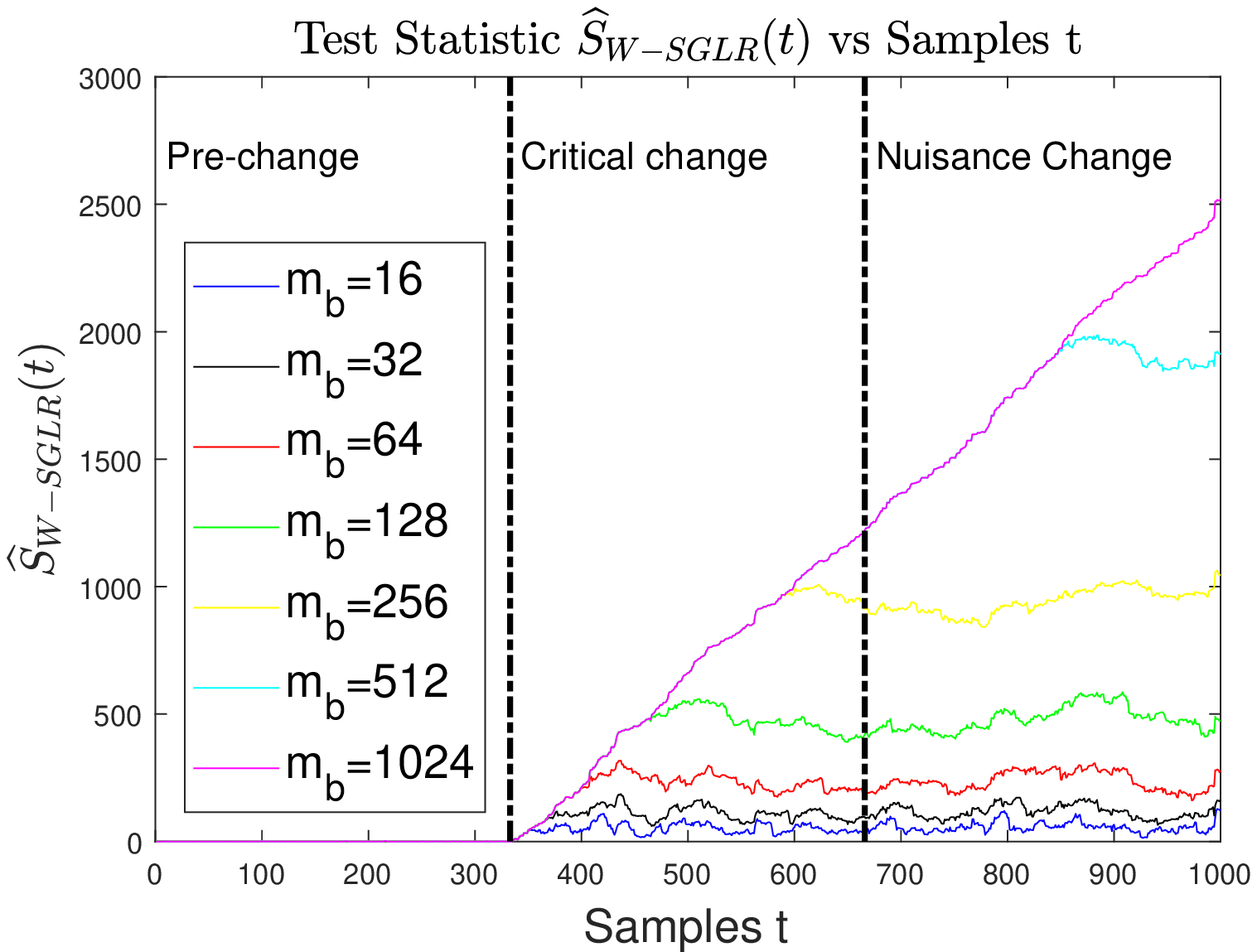}
			\caption{$\nu_c=333$, $\nu_n=666$}
			\label{fig:composite_nuisance_mean_nui666_crit333}
		\end{subfigure}
		\begin{subfigure}[b]{\columnwidth}
			\centering
			\includegraphics[width=0.85\columnwidth]{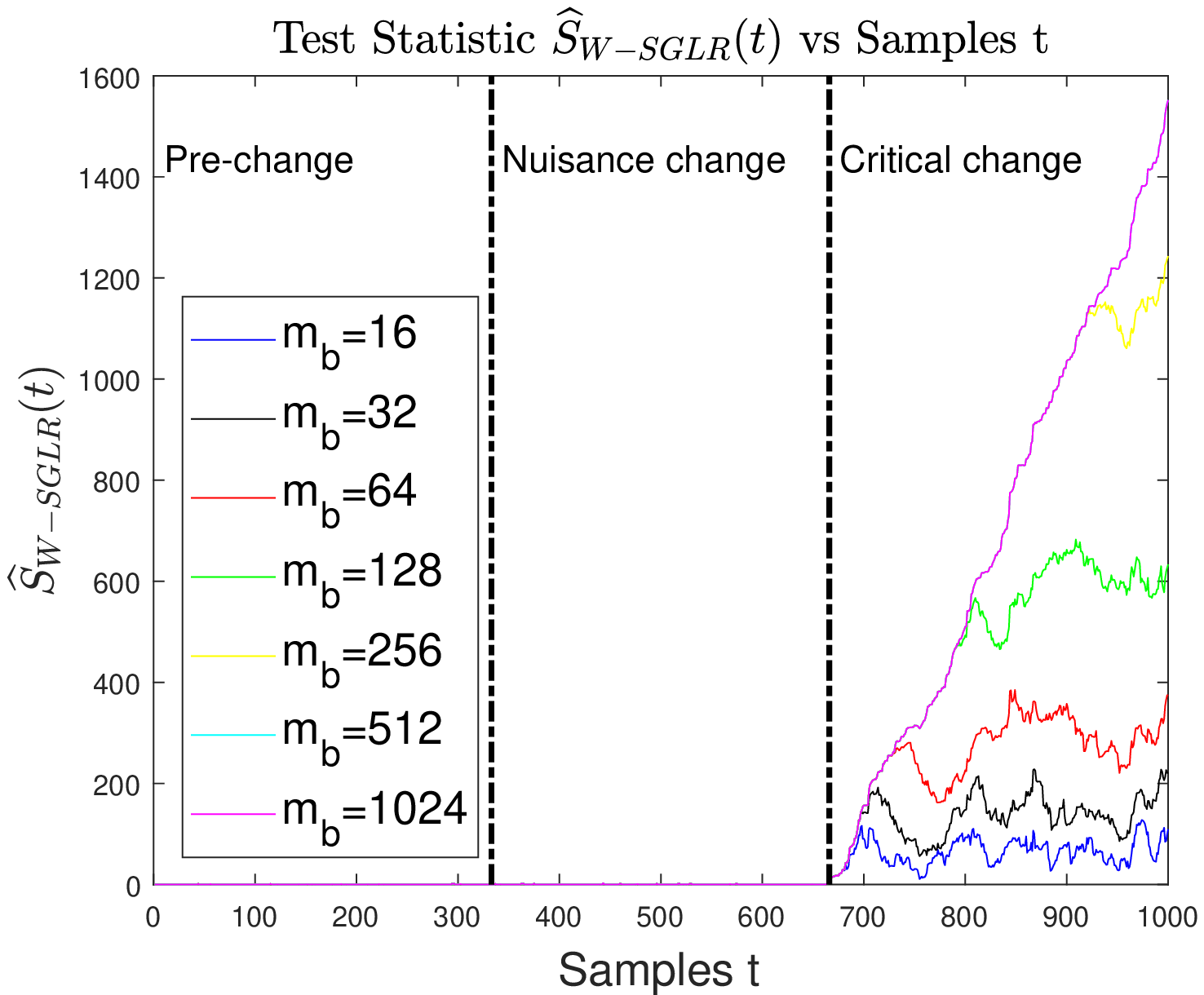}
			\caption{$\nu_c=666$, $\nu_n=333$}
			\label{fig:composite_nuisance_mean_nui333_crit666}
		\end{subfigure}
		\caption{The generalized W-SGLR test statistics $\widehat{S}_{\text{W-SGLR}}(t)$ with $f=\calN(0,1)$, $f_n=\calN(2,1)$, $g=\calN(0,10)$, $g_n=\calN(2,10)$ and $I=3.34$.}
		\label{fig:composite_nuisance_mean}
	\end{figure}
	
	In this set of simulations, we let $f=\calN(0,1)$, the critical change to be a change in variance where $g=\calN(0,\theta^2)$, and the nuisance change to be a change in the mean where $f_n=\calN(2,1)$ and $g_n=\calN(2,\theta_n^2)$. The parameters $\theta=\theta_n=10$ are unknown to the change detection algorithms. This corresponds to the case where the transmission power is unknown in the problem of spectrum sensing\cite{poor08qcdspectrum}. We ran the simulations with two change-point configurations to demonstrate the behavior of the generalized W-SGLR test statistic used in $\widehat{\tau}_{\text{W-SGLR}}$ as described in \eqref{eqn:modified_SGLR1} and \eqref{eqn:modified_SGLR2} for window-sizes $m_b=16,32,64,128,256,512$ and $1024$. In \figref{fig:composite_nuisance_mean_nui666_crit333}, we set the critical change point to be $\nu_{c}=333$ and nuisance change point to be $\nu_n=666$. It can be observed that our proposed generalized W-SGLR test statistic remains low during the pre-change regime, increases in the post-change regime and continue to increase in the nuisance post-change regime when the window is sufficiently large. This demonstrates that our stopping time $\widehat{\tau}_{\text{W-SGLR}}$ is able to detect the critical change even in the nuisance critical change region. In \figref{fig:composite_nuisance_mean_nui333_crit666}, we set the critical change point to be $\nu_{c}=666$ and nuisance change point to be $\nu_n=333$. %
	
	We see that our stopping time $\widehat{\tau}_\text{W-SGLR}$ is effective in detecting critical changes while ignoring the nuisance change in pre-change regime for window sizes as small as $m_b=16$. In practice, we can use graphs like  \cref{fig:composite_nuisance_mean_nui666_crit333} and \cref{fig:composite_nuisance_mean_nui333_crit666} to compare if the increase in the test-statistic after the critical change is discernible from the test-statistic in the pre-change regime. This would provide assistance in determining if the choice window-size is suitable.

\begin{figure}[!htbp]
		\centering
		\includegraphics[width=0.85\columnwidth]{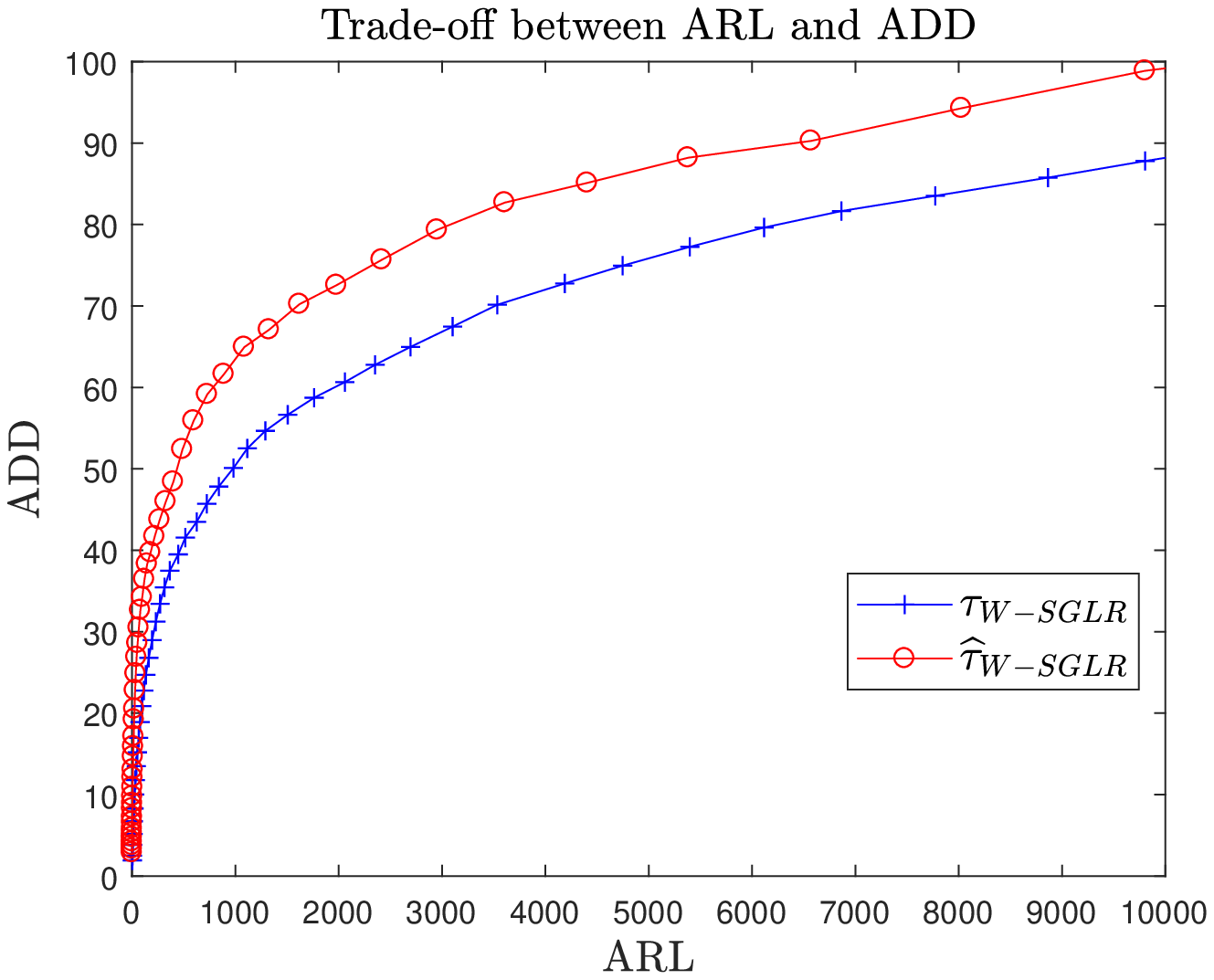}
		\caption{Comparison of trade-off performance for the proposed stopping time $\widehat{\tau}_{\text{W-SGLR}}$ and $\tau_{\text{W-SGLR}}$.}
		\label{fig:comparisonwithcomposite}
	\end{figure}

	Next, we compare the generalized W-SGLR stopping time with the W-SGLR stopping time. In our simulations, our signal is generated using the following distributions $f=\mathcal{N}(0,1)$, $f_n=\mathcal{N}(0,2)$, $g=\mathcal{N}(\theta,1)$, $g_n=\mathcal{N}(\theta_n,2)$. Here we set $\theta=\theta_n=2$ and assume that the condition that $\widehat{\theta}\in\text{Int}(\Theta)$ is always satisfied. We generate a signal of length $2^{16}=65,536$ and independently select the nuisance change point and critical change point with uniform probability on the $2^{16}$ possible data points. A total of $2^{12}=4096$ signals are generated. We compare the trade-off between the \gls{ADD} and the \gls{ARL} of the proposed W-SGLR stopping time when $\theta$ and $\theta_n$ are known against the generalized W-SGLR stopping time when $\theta$ and $\theta_n$ are unknown in \figref{fig:comparisonwithcomposite}. We observe that the generalized W-SGLR stopping time has a higher \gls{ADD} as compared to the W-SGLR stopping time. Our experiments suggest that the difference in \gls{ADD} is bounded as the \gls{ARL} becomes large.

	\subsection{Real Data}\label{subsec:RealData}
	In this subsection, we test our proposed stopping time $\tau_{\text{W-SGLR}}$ on the Case Western Reserve University Bearing Dataset \cite{smith2015rolling}. The dataset is collected from experiments conducted using an electric motor with accelerometer data measured at locations near to and remote from the motor bearings. Samples were collected at 12 KHz. We pre-process the signal by de-trending the signal using a first order finite difference: for each signal sample time $t$, let
	$
	X_t = Y_{t} - Y_{t-1},
	$
	where $Y_t$ is the observed raw signal sample at time $t$.
	
	\begin{figure}[!htb]
		\centering
		\includegraphics[width=0.85\columnwidth]{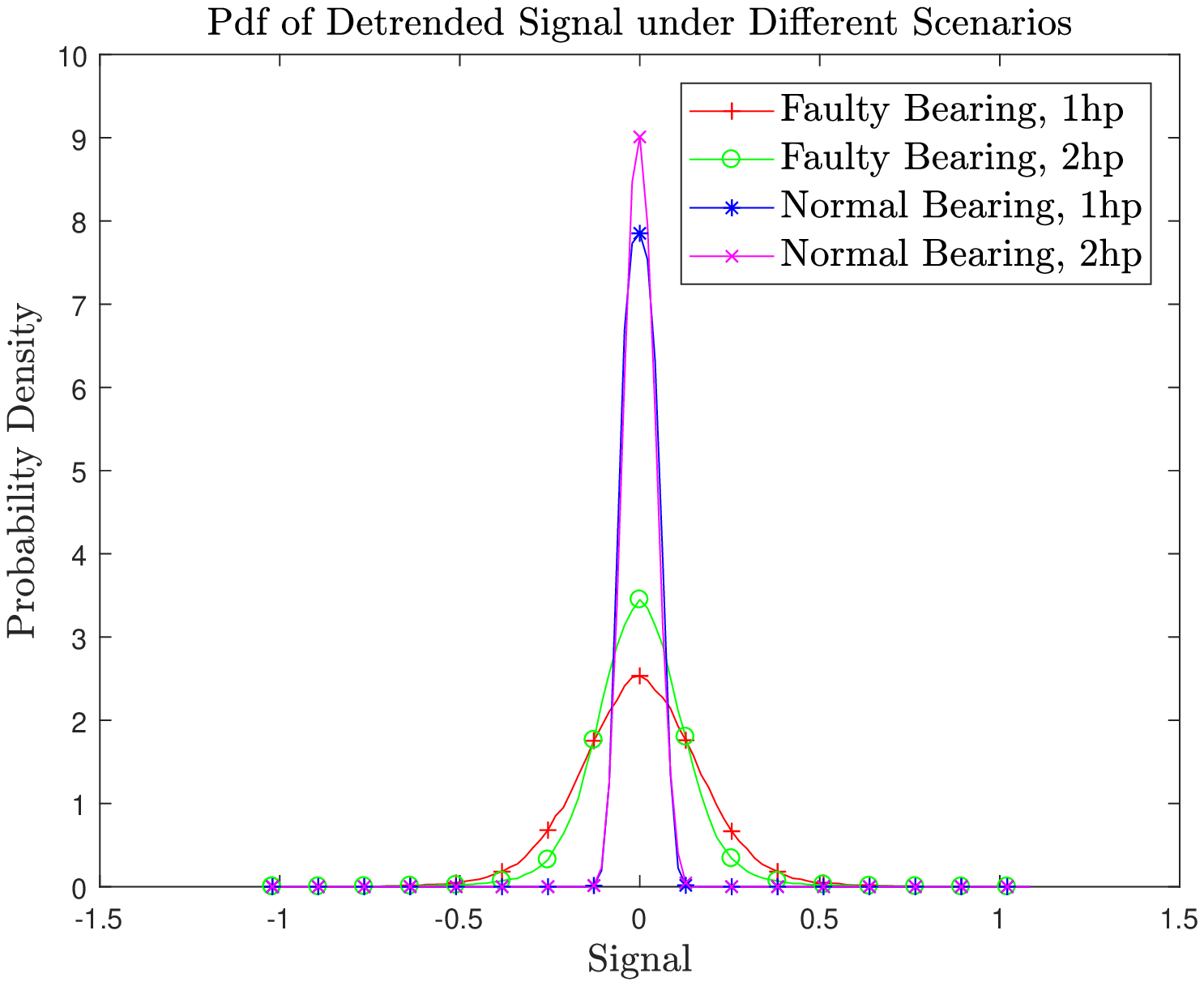}
		\caption{Plots of the pdfs of the observed de-trended signal in different scenarios.}
		\label{fig:pdf_of_different_scenario}
	\end{figure}
	
	We consider signals $X_t$ obtained at a motor load of 1hp and 2hp with normal bearings and also faulty bearings with a 0.007-inch fault diameter. We assume that the critical change would be the transition from a normal to faulty bearing, and a nuisance change would be a change in the motor load. We use the first 12,000 samples as training data to build a model for each of the following scenarios: normal bearings under a motor load of 1hp, normal bearings under a motor load of 2hp, faulty bearings under a motor load of 1hp, and faulty bearings under a motor load of 2hp. \figref{fig:pdf_of_different_scenario} shows the learned distributions of the de-trended signals observed in each scenario.
	
	\begin{figure}[!htb]
		\begin{subfigure}{\columnwidth}
			\centering
			\includegraphics[width=0.85\columnwidth]{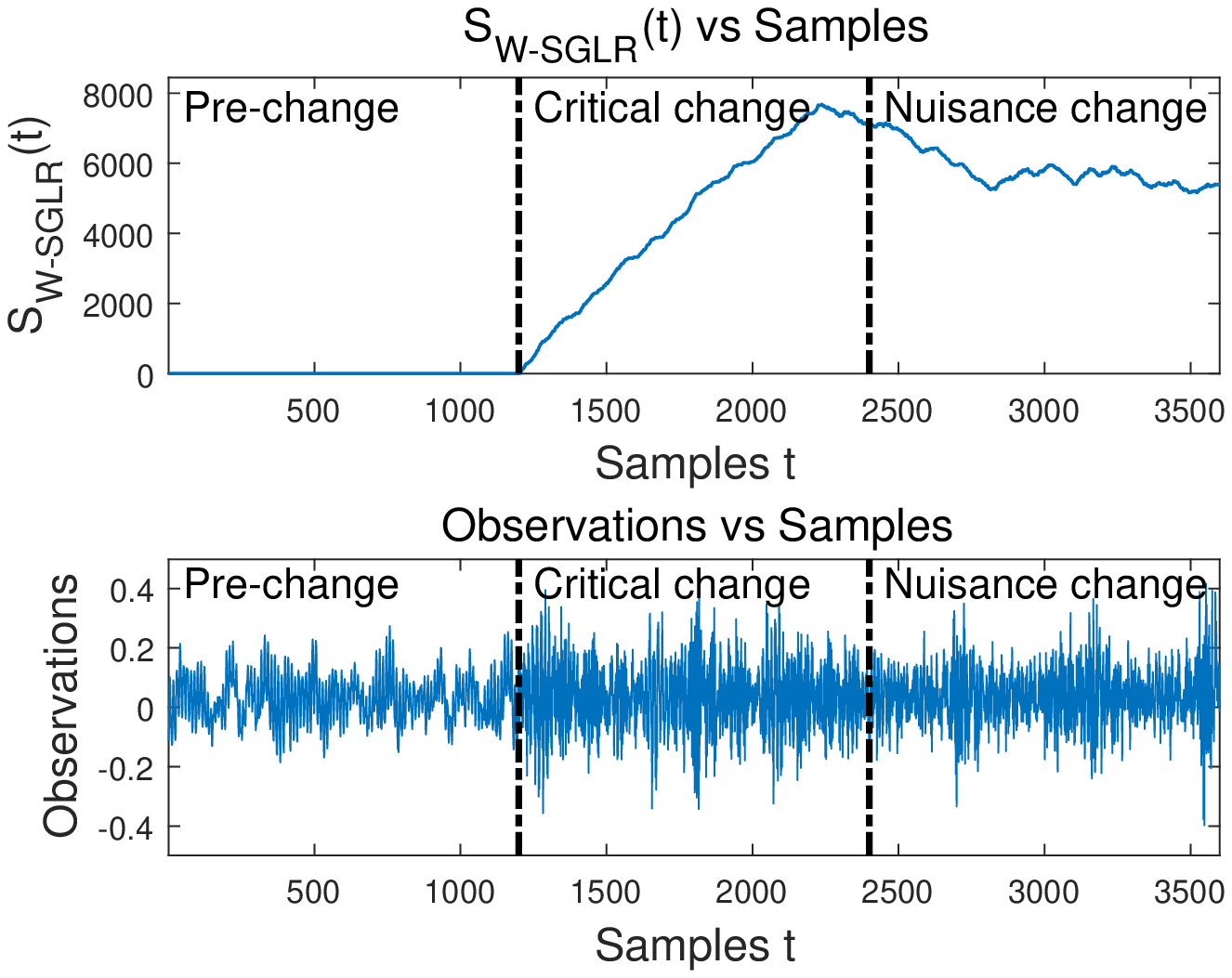}
			\caption{$\nu_c=1200$, $\nu_n=2400$}
			\label{fig:real_example_cri_nui}
		\end{subfigure}
		\begin{subfigure}{\columnwidth}
			\centering
			\includegraphics[width=0.85\columnwidth]{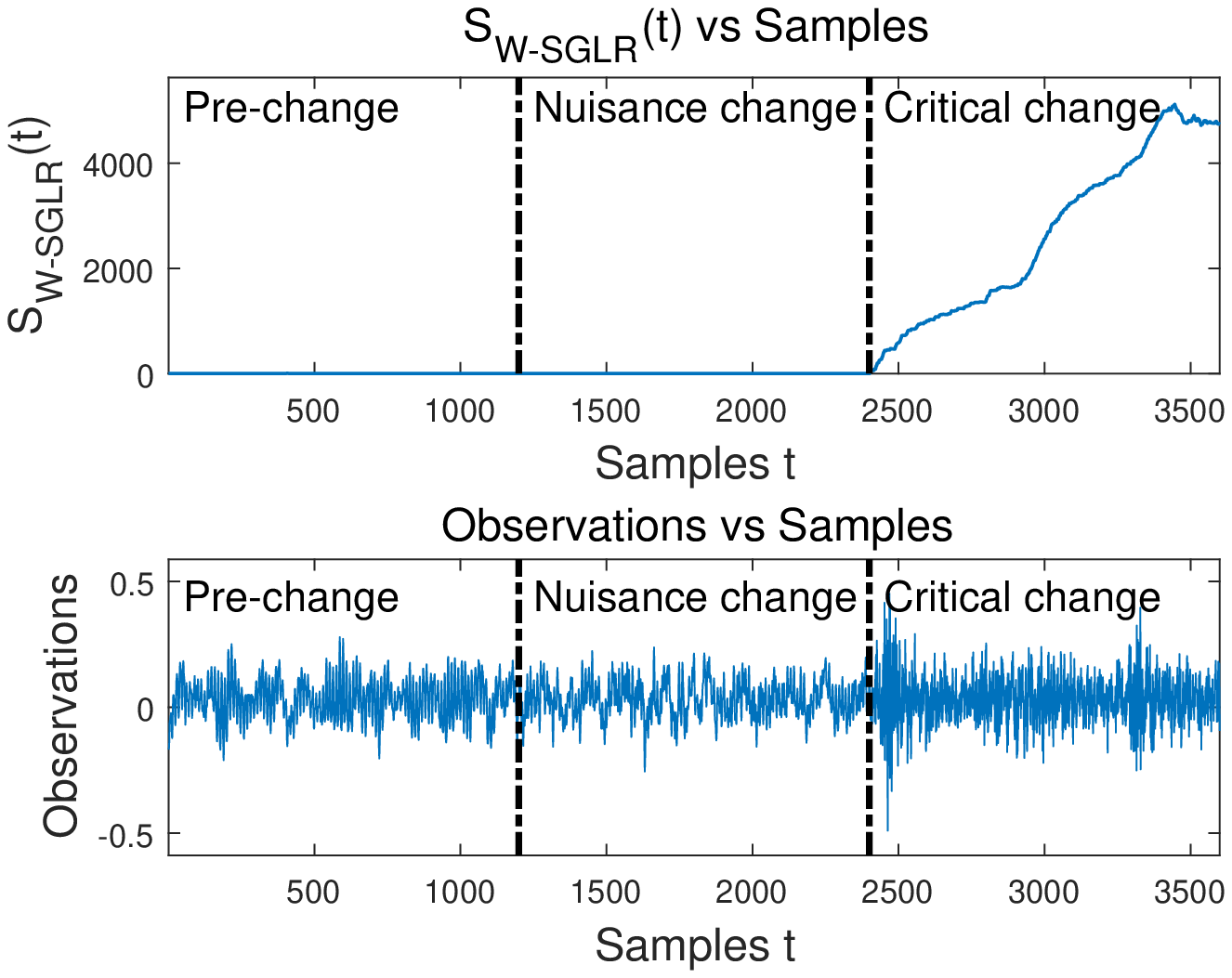}
			\caption{$\nu_c=2400$, $\nu_n=1200$}
			\label{fig:real_example_nui_cri}
		\end{subfigure}
		\caption{Examples of the W-SGLR test statistic with $m_b=1024$.}\label{fig:real_example}
	\end{figure}

	There are two challenges faced in testing our proposed stopping time on real data: (i) we lack theoretical results for the \gls{ARL} of the 2-stage stopping times for the selection of appropriate thresholds for comparison and (ii) real run-to-failure data is difficult to obtain. We divide the remaining samples into 3 disjoint sets to address the above challenges. 
	
	For the first set, we create a training set of 1000 signals each with length 36,000 with a randomly selected nuisance change point $\nu_n$ for each signal such that there is a period of $\nu_n-1$ samples for a normal bearing under a motor load of 1hp, and a period of $36,000-\nu_n+1$ samples for a normal bearing under a motor load of 2hp. We select appropriate thresholds for each of the stopping times so that the empirical ARL varies between 1200 and 18,000. 
	
	The next two sets are testing sets. We create 1000 signals of length 3600 each with (i) a period of 1200 samples for a normal bearing under a motor load of 1hp, which transitions to (ii) a period of 1200 samples for a normal bearing under a motor load of 2hp, which finally transitions to (iii) a period of 1200 samples for a faulty bearing under a motor load of 2hp. 
	
	Similarly, we create 1000 signals for the scenario where a normal bearing under a motor load of 1hp transitions to a faulty bearing under a motor load of 1hp and finally a faulty bearing under a motor load of 2hp.

	Finally, we apply the selected thresholds obtained from the first training set to the two testing sets to compute the stopping times' empirical ADD performance. The window size of $m_b=8$ for the FMA stopping time is selected to minimize its empirical ADD on the test set. 
	For this dataset, if $m_b$ is chosen to be $128$ or $1024$, the empirical ADD of the FMA stopping time becomes much larger compared to the empirical ADD of the W-SGLR and 2-stage stopping times. Thus, we only present the ARL-ADD trade-off for $m_b=8$.
	
	In \cref{fig:real_example_cri_nui,fig:real_example_nui_cri}, we present some examples of the performance of the W-SGLR test statistic. It can be seen that in both cases, the test-statistic remains low before the bearing failure and quickly rises after the bearing fails even as the motor load changes. 
	
\begin{figure}[!htb]
		\begin{subfigure}[b]{\columnwidth}
			\centering
			\includegraphics[width=0.85\columnwidth]{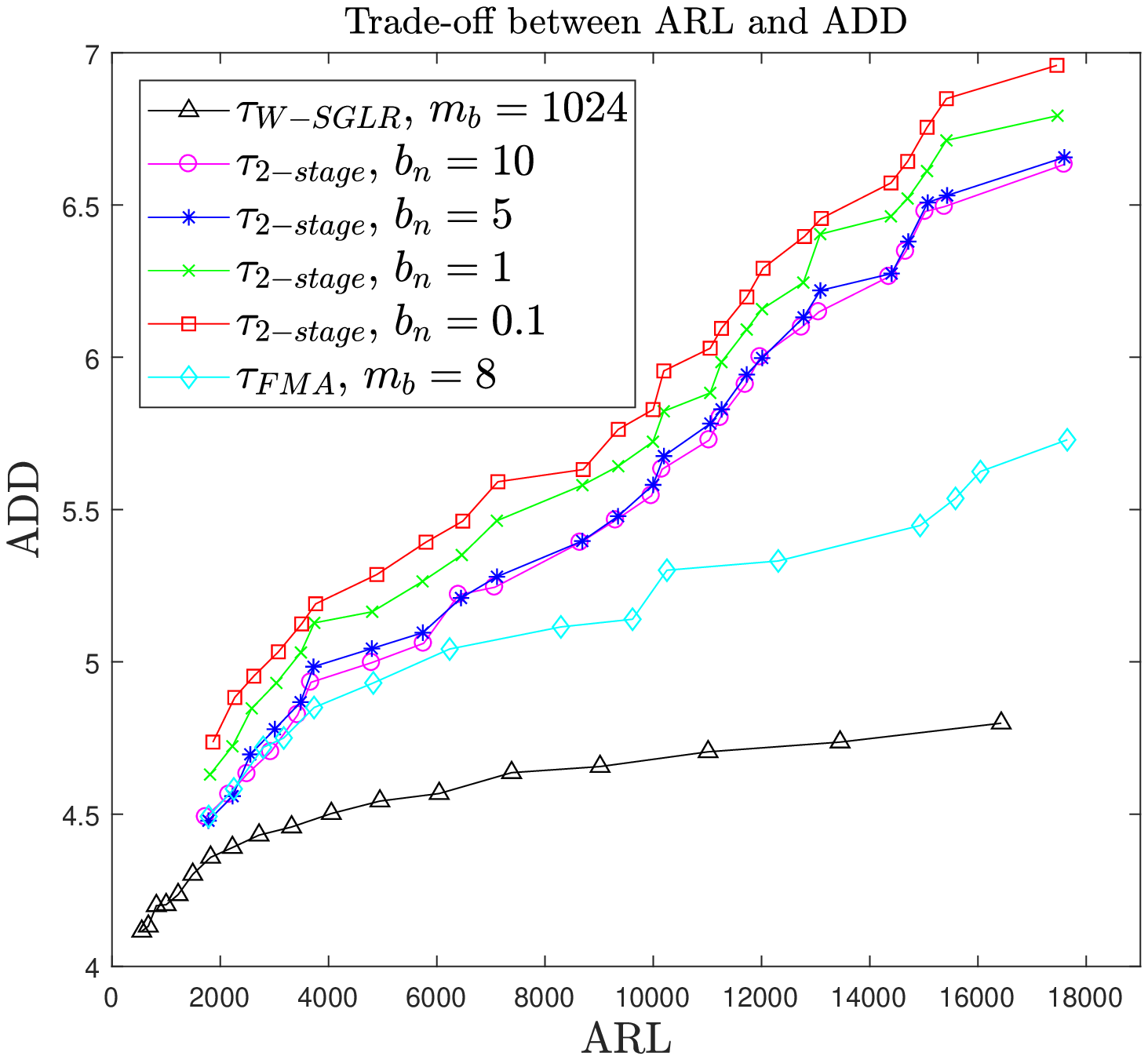}
			\caption{$\nu_n=1200$, $\nu_c=2400$}
			\label{fig:real_example_add_vs_arl_vs_2stage_nui_crit}
		\end{subfigure}
		\ \\ 
		\begin{subfigure}[b]{\columnwidth}
			\centering
			\includegraphics[width=0.85\columnwidth]{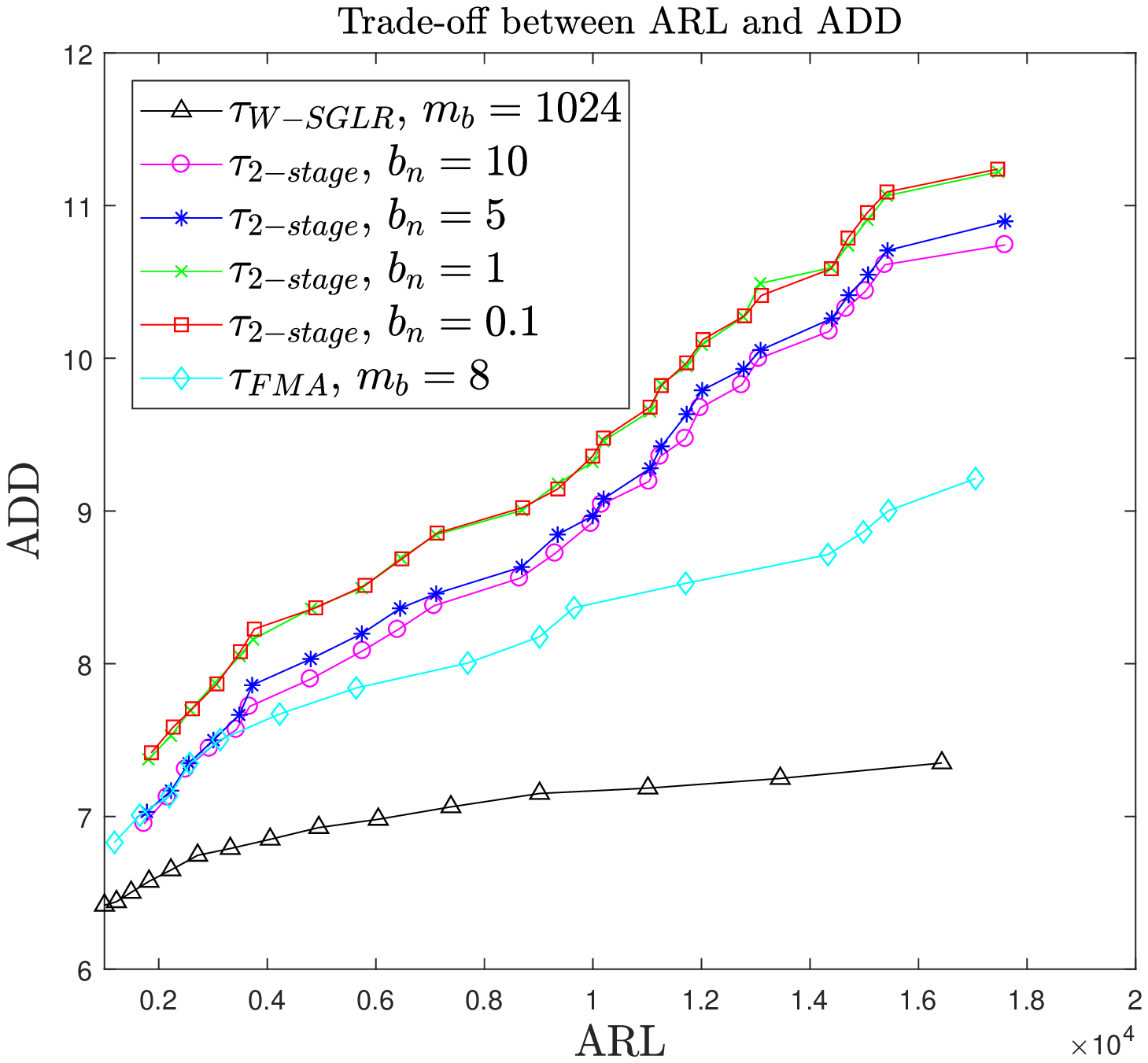}
			\caption{$\nu_n=2400$, $\nu_c=1200$}
			\label{fig:real_example_add_vs_arl_vs_2stage_crit_nui}
		\end{subfigure}
		\caption{Trade-off between the empirical ADD and ARL.}
		\label{fig:real_example_add_vs_arl_vs_2stage}
	\end{figure}

	In \cref{fig:real_example_add_vs_arl_vs_2stage_nui_crit} and \cref{fig:real_example_add_vs_arl_vs_2stage_crit_nui}, we present the trade-off between the empirical ADD and ARL for the proposed W-SGLR stopping time with $m_b=1024$, the 2-stage stopping times  with different thresholds $b_n$ and the FMA stopping time. It can be seen that our proposed stopping time $\tau_{\text{W-SGLR}}$ achieves better ADD-ARL trade-off compared to the other stopping times. 
	However, as the KL divergences $\KLD{g_n}{f_n},\KLD{g}{f_n},\KLD{g_n}{f},\KLD{g}{f}$ are large, the empirical ADD for all the algorithms remains low across the range of ADD tested. In this case, the reduction in empirical ADD is small, between $1$ to $4$ samples, over the range of ARLs tested. In terms of computational complexity, up till sample $t$, the W-SGLR stopping time requires $O(m_bt)$ operations\cite{lau2019quickest}, which is slightly more than both the two-stage stopping time and the FMA stopping time, both of which require $O(t)$ operations. Thus, for applications that have limited computational resources and large differences in their pre and post-change distributions, we may want to consider using the FMA or the 2-stage stopping time as the degradation in performance is small.
	
	\section{Discussions and Conclusions}\label{sec:conclusion}
	We have studied the non-Bayesian QCD problem where the signal may be subjected to a nuisance change. We proposed the W-SGLR stopping time that quickly detects the critical change while ignoring the nuisance change. The limited window size ensures that the W-SGLR stopping time does not require increasing computational resources as more samples are observed. We also derived the stopping time's asymptotic behavior and showed that it is asymptotically optimal under mild technical assumptions. A generalized W-SGLR stopping time is also proposed for the case where the critical and nuisance post-change distributions are unknown but belong to a parametrized family. Numerical simulations and experiments on a real dataset demonstrated that the W-SGLR stopping time achieves better ADD-ARL trade-off than various other competing stopping times.

	In this paper, we have assumed that if both the critical and nuisance changes occur, the eventual distribution that generates the signal is the same, regardless of which change comes first. A more general model would be to allow the eventual distribution to depend on the order of the change points. An easy generalization of the W-SGLR stopping would be to include all the different eventual distributions into the numerator of \cref{eqn:SGLR_teststats}. The asymptotic trade-off between the WADD and ARL can be derived using similar techniques in \cref{subsec:properties_stats}. However, deriving the conditions for asymptotic optimality of this stopping time is more complicated and would be a possible direction for future research. 
	
	Another possible future research direction is to consider a modification of the W-SGLR stopping time for the TCD problem under the possibility of a nuisance change. As the performance metrics of the TCD problem are different from the QCD problem, its asymptotic trade-off between the worst-case false alarms and missed detection within a specified window needs to be studied. 
	Also, as the FMA performs well in the TCD problem, it will be interesting to consider if the FMA stopping time can be adapted to solve our QCD problem.
	
	\appendices
	
	\section{Proof of Propositions \ref{prop:converge_to_post_nui_LR} and \ref{prop:converge_to_pre_nui_LR}}\label[appendix]{sec:AppProp1}
	
	We start off with some notation definitions. Let $L_i=\log\tfrac{f_n(X_i)}{f(X_i)}$. For any $N\geq0$, let
	$
	L_{i,>N}=L_i\mathbf{1}_{\{|L_i|>N\}}\ \text{and}\  L_{i,\leq N}=L_i\mathbf{1}_{\{|L_i|\leq N\}}.
	$
	For any $k,t\in\mathbb{N}$ such that $k\leq t$, we define the following averages:
	\begin{align}
	\begin{split}&\overline{L^{k:t}}=\tfrac{1}{t-k+1}\sum_{i=k}^tL_i,\ \overline{L^{k:t}_{>N}}=\tfrac{1}{t-k+1}\sum_{i=k}^tL_{i,>N},\\ 
	&\quad\quad\quad\overline{L^{k:t}_{\leq N}}=\tfrac{1}{t-k+1}\sum_{i=k}^tL_{i,\leq N}\end{split}
	\end{align}
	We have
	\begin{align}\label{eqn:mean_truncation}
	\overline{L^{k:t}}=\overline{L^{k:t}_{>N}}+\overline{L^{k:t}_{\leq N}}.
	\end{align}
	For the case where $k>t$, we let $\overline{L^{k:t}}=\overline{L^{k:t}_{>N}}=\overline{L^{k:t}_{\leq N}}=0$. Finally, we define the random variable
	\begin{align}
	V_{k,t}=\arg\max_{k\leq j\leq t+1}\prod_{i=k}^{j-1}f(X_i)\prod_{i=j}^{t}f_n(X_i).
	\end{align}
	
	An outline of the proof of Propositions \ref{prop:converge_to_post_nui_LR} and \ref{prop:converge_to_pre_nui_LR} is as follows. Lemma~\ref{lem:no_nui_change_gap_small} and Lemma~\ref{lem:truncation_error_bound} provide the results required for controlling the error bound in Proposition~\ref{prop:converge_to_post_nui_LR}. Similarly, Lemma~\ref{lem:nui_change_gap_small} and Lemma~\ref{lem:truncation_error_bound} provide the results required for controlling the error bound in Proposition~\ref{prop:converge_to_pre_nui_LR}. The Lemmas~\ref{lem:no_nui_change_gap_small},~\ref{lem:nui_change_gap_small} and~\ref{lem:truncation_error_bound} require that decay in the tail probabilities of the average log-likelihood ratio $\tfrac{1}{n}\sum_{i=k}^{k+n-1}\log\tfrac{f_n(X_i)}{f(X_i)}$ to be at most $O(n^{-2})$, which is shown in \cref{lem:fourth_order_tail_bound}.
	
	\begin{Lemma_A}\label{lem:fourth_order_tail_bound}
		For any $\nu_c,\nu_n,k,l,n\in\mathbb{N}$ such that $\nu_c\leq k$ and $\max\{\nu_n,\nu_c\}\leq l$, and $\epsilon>0$, we have
		\begin{align}
		\P{\infty,\nu_c}(\left|\overline{L^{k:k+n-1}}-\rho_{g}\right|\geq \epsilon)  & \leq \tfrac{K_{g}}{\epsilon^4 n^2}, \label{4order1}   \\
		\P{\nu_n,\nu_c}(\left|\overline{L^{l:l+n-1}}-\rho_{g_n}\right|\geq \epsilon) & \leq \tfrac{K_{g_n}}{\epsilon^4 n^2}, \label{4order2}
		\end{align}
		where $K_g=\omega_g^4+\tfrac{3}{2}\sigma_g^4$ and $K_{g_n}=\omega_{g_n}^4+\tfrac{3}{2}\sigma_{g_n}^4$.
	\end{Lemma_A}
	\begin{IEEEproof}
		As the proof is elementary, we omit it here and refer the reader to the extended version in \cite{lau2019quickestarxiv}.	
	\end{IEEEproof}

	From \cref{lem:fourth_order_tail_bound}, for any $\nu_c \leq k \leq v \leq t < \infty$, we have for $\rho_g < 0$,
	\begin{align}
	& \P{\infty,\nu_c}(V_{k,t}=v) \nn
	& \leq\P{\infty,\nu_c}(\left|\overline{L^{v:t}}-\rho_g\right|\geq |\rho_g|)\nn
	& \leq\tfrac{K_{g}}{|\rho_g|^4 (t-v+1)^2}.\label{Vkt_bound}
	\end{align}
	Similarly, for $\max\{\nu_c,\nu_n\} \leq k < v \leq t+1 <\infty$, and $\rho_{g_n}>0$, we have
	\begin{align}
	\P{\nu_n,\nu_c}(V_{k,t}=v) & \leq\tfrac{K_{g}}{|\rho_{g_n}|^4 (v-k)^2}.\label{Vktn_bound}
	\end{align}
	
	For the next two lemmas, we use bounds on the tail probability of $\overline{L^{k:t}_{>N_{g}}}$ to derive asymptotic properties of the random variable $V_{k,t}$ under the distributions $\P{\nu_n,\nu_c}$ and $\P{\infty,\nu_c}$ for any $\nu_n,\nu_c\in\mathbb{N}$.
	
	\begin{Lemma_A}\label{lem:no_nui_change_gap_small}
		For any $0<c< 1$, $\nu_c\leq k<\infty$, and $\rho_g < 0$,  we have
		\begin{align}\label{eqn:no_nui_change_gap_small}
		\lim_{t\to\infty}\P{\infty,\nu_c}(\tfrac{t-V_{k,t}+1}{t-k+1}> c)=0.
		\end{align}
	\end{Lemma_A}
	\begin{IEEEproof}
		We have
		\begin{align}
		& \P{\infty,\nu_c}(\tfrac{t-V_{k,t}+1}{t-k+1}> c)\nn
		& \leq\sum_{v=k}^{\floor{(1-c)(t+1)+ck}}\P{\infty,\nu_c}(V_{k,t}=v)\nn
		& \leq\sum_{v=k}^{\floor{(1-c)(t+1)+ck}}\tfrac{K_{g}}{|\rho_g|^4 (t-v+1)^2}\label{eqn:use_fourth_moment_tail_bound} \\
		& \leq\left(\floor{(1-c)(t+1)+ck}-k+1\right)\tfrac{K_{g}}{|\rho_g|^4 (t-k+1)^2}\to 0, \nonumber
		\end{align}
		as $t\to\infty$. The inequality \eqref{eqn:use_fourth_moment_tail_bound} follows from \cref{Vkt_bound}. The proof is now complete.
	\end{IEEEproof}
	
	\begin{Lemma_A}\label{lem:nui_change_gap_small}
		{For any $0<c< 1$ and $\nu_c\leq k<\infty$, let  $k'=\max\{k,\nu_n\}$. If $\rho_{g_n}>0$, we have
			\begin{align}\label{eqn:nui_change_gap_small}
			{\lim_{t\to\infty}\P{\nu_n,\nu_c}(\tfrac{V_{k,t}-k'}{t-k+1}> c)=0}.
			\end{align}
		}
	\end{Lemma_A}
	\begin{IEEEproof}
		For $t>\tfrac{k-k'}{c}+k-1$, we have $c(t-k+1)+k'>k$ and
		\begin{align}
		& \P{\nu_n,\nu_c}(\tfrac{V_{k,t}-k'}{t-k+1}> c)\nn
		& =\sum_{v=\ceil{c(t-k+1)+k'}}^{t+1}\P{\nu_n,\nu_c}(V_{k,t}=v)\nn
		& \leq\sum_{v=\ceil{c(t-k+1)+k'}}^{t+1}\tfrac{K_{g_n}}{|\rho_{g_n}|^4 (v-k)^2}\label{eqn:nui_use_fourth_moment_sum_bound} \\
		& \leq\left(t+1-\ceil{c(t-k+1)+k'}+1\right)\tfrac{K_{g_n}}{|\rho_{g_n}|^4 (\ceil{c(t-k+1)}+k'-k)^2}\to 0, \nonumber
		\end{align}
		as $t\to\infty$, where \cref{eqn:nui_use_fourth_moment_sum_bound} follows from \cref{Vktn_bound}. The proof is now complete.
	\end{IEEEproof}
	\begin{Lemma_A}\label{lem:truncation_error_bound}
		Suppose that $\rho_{g}<0$ and $\rho_{g_n}>0$. For any $\epsilon,\delta>0$ and $k\geq \max\{\nu_n,\nu_c\}$, there exist $N_{g_n},N_g\in\mathbb{N}$ such that for any $t\geq k$,
		\begin{align}
		\P{\infty,\nu_c}(\left|\overline{L^{V_{k,t}:t}_{>N_{g}}}\right|\geq \epsilon)    & \leq\delta,\label{eqn:truncate_error_1} \\
		\P{\nu_n,\nu_c}(\left|\overline{L^{k:V_{k,t}-1}_{>N_{g_n}}}\right|\geq \epsilon) & \leq\delta.\label{eqn:truncate_error_2}
		\end{align}
		
	\end{Lemma_A}
	\begin{IEEEproof}
		Given any $\epsilon,\delta>0$, since the fourth moment of $\log\tfrac{f_n(X)}{f(X)}$ exists, by the monotone convergence theorem, there exists $N_g$ and $N_{g_n}$ such that
		\begin{align}
		\E{g_n}[\left|\log \tfrac{f_n(X)}{f(X)}\right|^4\mathbf{1}_{\left\{\left|\log \tfrac{f_n(X)}{f(X)}\right|>N_{g_n}\right\}}]\leq\tfrac{\epsilon^4\delta^4}{M^4},\label{LNgn} \\
		\E{g}[\left|\log \tfrac{f_n(X)}{f(X)}\right|^4\mathbf{1}_{\left\{\left|\log \tfrac{f_n(X)}{f(X)}\right|>N_{g}\right\}}]\leq\tfrac{\epsilon^4\delta^4}{M^4},
		\end{align}
		where $M=\sum_{v=1}^{\infty}\left(\tfrac{K_{g}}{|\rho_g|^4 v^2}\right)^{\tfrac{3}{4}} < \infty$. Applying Markov's inequality, we obtain
		\begin{align}\label{eqn:tallprobability_control}
		\P{\infty,\nu_c}(\left|\overline{L^{V_{k,t}:t}_{>N_{g}}}\right|\geq \epsilon)\leq\ofrac{\epsilon}\E{\infty,\nu_c}[\left|\overline{L^{V_{k,t}:t}_{>N_{g}}}\right|].
		\end{align}
		
		Next, we derive an upper bound for $\E{\infty,\nu_c}[\left|\overline{L^{V_{k,t}:t}_{>N_{g}}}\right|]$. For any $v \leq t$, we have
		\begin{align}
		\E{\infty,\nu_c}[\left|\overline{L^{v:t}_{>N_{g}}}\right|^4] & =\E{\infty,\nu_c}[\left|\tfrac{1}{t-v+1}\sum_{i=v}^t L_{i,> N_g}\right|^4]\nn
		& \leq\tfrac{1}{t-v+1}\sum_{i=v}^t\E{\infty,\nu_c}[\left|L_{i,>N_g}\right|^4]\label{jensen} \\
		& \leq \tfrac{\epsilon^4\delta^4}{M^4},\label{bddLvt}
		\end{align}
		where \cref{jensen} follows from Jensen's inequality, and \cref{bddLvt} follows from \cref{LNgn}. We obtain
		\begin{align}
		\E{\infty,\nu_c}[\left|\overline{L^{V_{k,t}:t}_{>N_{g}}}\right|] & =\sum_{v=k}^{t+1}\E{\infty,\nu_c}[\left|\overline{L^{V_{k,t}:t}_{>N_{g}}}\right|\indicator{V_{k,t}=v}]\nn
		& =\sum_{v=k}^{t}\E{\infty,\nu_c}[\left|\overline{L^{v:t}_{>N_{g}}}\right|\indicator{V_{k,t}=v}]\label{eqn:reduced_sum_1}                           \\
		& \leq\sum_{v=k}^{t}\E{\infty,\nu_c}[\left|\overline{L^{v:t}_{>N_{g}}}\right|^4]^{\ofrac{4}}\P{\infty,\nu_c}(V_{k,t}=v)^{\tfrac{3}{4}}\label{holder} \\
		& \leq\tfrac{\epsilon\delta}{M}\sum_{v=k}^{t}(\P{\infty,\nu_c}(V_{k,t}=v))^{\tfrac{3}{4}}, \label{bdd1}                                               \\
		& \leq \tfrac{\epsilon\delta}{M}\sum_{v=k}^{t}\left(\tfrac{K_{g}}{|\rho_g|^4 (t-v+1)^2}\right)^{\tfrac{3}{4}}\label{bdd2}                              \\
		& \leq \epsilon\delta. \label{ELNg_bdd}
		\end{align}
		where \eqref{eqn:reduced_sum_1} is because $\overline{L^{t+1:t}_{>N_{g}}}=0$, \cref{holder} follows from H\"older's inequality, \cref{bdd1} from \cref{bddLvt}, \cref{bdd2} from \cref{Vkt_bound}, and \cref{ELNg_bdd} from the definition of $M$. From \eqref{eqn:tallprobability_control}, we have
		$
		\P{\infty,\nu_c}(\left|\overline{L^{V_{k,t}:t}_{>N_{g}}}\right|\geq \epsilon)\leq \tfrac{\epsilon\delta}{\epsilon} = \delta,
		$
		and \cref{eqn:truncate_error_1} is proved. The proof of \eqref{eqn:truncate_error_2} is similar and the lemma is proved.
	\end{IEEEproof}

	\subsection{Proof of Proposition \ref{prop:converge_to_post_nui_LR}}
	It suffices to show that for any $\epsilon,\delta>0$, there exists $T$ such that for all $t\geq T$ we have
	\begin{align}\label{eqn:prop_3_1}
	\P{\infty,\nu_c}(\left|\tfrac{\log \Lambda(k,t)}{t-k+1}-\tfrac{1}{t-k+1}\sum_{i=k}^t\log \tfrac{g(X_i)}{f(X_i)}\right|\geq\epsilon)\leq\delta.
	\end{align}
	For any $N\geq0$ and $c>0$, the left-hand side of \eqref{eqn:prop_3_1} becomes
	\begin{align}
	& \P{\infty,\nu_c}(\tfrac{1}{t-k+1}\left|\sum_{i=V_{k,t}}^t\log L_i \right|\geq\epsilon)\nn
	\begin{split}& \leq \P{\infty,\nu_c}(\tfrac{t-V_{k,t}+1}{t-k+1}\left| \overline{L^{V_{k,t}:t}_{>N}}\right|\geq\tfrac{\epsilon}{2})\\&\quad+ \P{\infty,\nu_c}(\tfrac{t-V_{k,t}+1}{t-k+1}\left| \overline{L^{V_{k,t}:t}_{\leq N}}\right|\geq\tfrac{\epsilon}{2})\end{split}\nn
	\begin{split}& \leq \P{\infty,\nu_c}(\left| \overline{L^{V_{k,t}:t}_{>N}}\right|\geq\tfrac{\epsilon}{2})\\&\quad+ \P{\infty,\nu_c}(\left\{\tfrac{t-V_{k,t}+1}{t-k+1}\left| \overline{L^{V_{k,t}:t}_{\leq N}}\right|\geq\tfrac{\epsilon}{2}\right\}\bigcap\left\{\tfrac{t-V_{k,t}+1}{t-k+1} \leq c\right\})\end{split}\nn
	& + \P{\infty,\nu_c}(\tfrac{t-V_{k,t}+1}{t-k+1}> c). \label{eqn:prop_3_1_bdd}
	\end{align}
	From \cref{lem:truncation_error_bound}, there exists $N$ such that $\P{\infty,\nu_c}(\left| \overline{L^{V_{k,t}:t}_{>N}}\right|\geq\tfrac{\epsilon}{2})\geq\tfrac{\delta}{2}$. Next, by choosing $c=\tfrac{\epsilon}{4N}$, we have
	$
	\P{\infty,\nu_c}(\left\{\tfrac{t-V_{k,t}+1}{t-k+1}\left| \overline{L^{V_{k,t}:t}_{\leq N}}\right|\geq\tfrac{\epsilon}{2}\right\}\bigcap\left\{\tfrac{t-V_{k,t}+1}{t-k+1} \leq c\right\})=0.
	$
	Finally, from Lemma~\ref{lem:no_nui_change_gap_small}, there exists $T$ such that for all $t\geq T$, we have
	$
	\P{\infty,\nu_c}(\tfrac{t-V_{k,t}+1}{t-k+1}> c)\leq\tfrac{\delta}{2}.
	$
	The right-hand side of \cref{eqn:prop_3_1_bdd} is then upper bounded by $\delta$, and the proof is complete.
	
	\subsection{Proof of Proposition \ref{prop:converge_to_pre_nui_LR}}
	It suffices to show that for any $\epsilon,\delta>0$, there exists $T$ such that for all $t\geq T$ we have
	\begin{align}\label{eqn:prop_4_1}
	\P{\nu_n,\nu_c}(\left|\tfrac{\log \Lambda(k,t)}{t-k+1}-\tfrac{1}{t-k+1}\sum_{i=k}^t\log
	\tfrac{g(X_i)}{f_n(X_i)}\right|\geq\epsilon)\leq\delta.\end{align}
	Let $k'=\max\{k,\nu_n\}$. The left-hand side of \eqref{eqn:prop_4_1} can be written as
	\begin{align}
	& \P{\nu_n,\nu_c}(\tfrac{1}{t-k+1}\left|\sum_{i=k}^{V_{k,t}-1}\log L_i\right|\geq\epsilon)\nonumber                                      \\
	& \leq\P{\nu_n,\nu_c}(\tfrac{1}{t-k+1}\left|\sum_{i=k}^{k'-1}\log L_i \right|\geq\tfrac{\epsilon}{2})\label{eqn:prop_4_pre_nui}           \\
	& \quad+\P{\nu_n,\nu_c}(\tfrac{1}{t-k+1}\left|\sum_{i=k'}^{V_{k,t}-1}\log L_i \right|\geq\tfrac{\epsilon}{2}).\label{eqn:prop_4_post_nui}\end{align}
	Applying Markov's inequality to \cref{eqn:prop_4_pre_nui}, there exists $T_1$ such that for all $t\geq T_1$, we have
	\begin{align}
	& \P{\nu_n,\nu_c}(\tfrac{1}{t-k+1}\left|\sum_{i=k}^{k'-1}\log L_i \right|\geq\tfrac{\epsilon}{2})\nonumber            \\
	& \leq \frac{2}{\epsilon(t-k+1)}\E{\nu_n,\nu_c}[\left|\sum_{i=k}^{k'-1}\log L_i \right|]
	<\tfrac{\delta}{3}.\nonumber
	\end{align}
	For any $N\geq 0$ and $c\geq 0$, \cref{eqn:prop_4_post_nui} becomes
	\begin{align}
	& \P{\nu_n,\nu_c}(\tfrac{1}{t-k+1}\left|\sum_{i=k'}^{V_{k,t}-1}\log L_i \right|\geq\tfrac{\epsilon}{2})\nonumber                                                                                                                                                                     \\
	\begin{split}&\leq  \P{\nu_n,\nu_c}(\tfrac{V_{k,t}-k'}{t-k+1}\left| \overline{L^{k':V_{k,t}-1}_{>N}}\right|\geq\tfrac{\epsilon}{4})\\&\quad+ \P{\nu_n,\nu_c}(\tfrac{V_{k,t}-k'}{t-k+1}\left| \overline{L^{k':V_{k,t}-1}_{\leq N}}\right|\geq\tfrac{\epsilon}{4})\end{split}\nonumber                                        \\
	\begin{split}&\leq  \P{\nu_n,\nu_c}(\left| \overline{L^{k':V_{k,t}}_{>N}}\right|\geq\tfrac{\epsilon}{4})\\&\quad+ \P{\nu_n,\nu_c}(\left\{\tfrac{V_{k,t}-k'}{t-k+1}\left| \overline{L^{V_{k,t}:t}_{\leq N}}\right|\geq\tfrac{\epsilon}{2}\right\}\bigcap\left\{\tfrac{V_{k,t}-k'}{t-k+1} \leq c\right\})\end{split}\nonumber \\
	& + \P{\nu_n,\nu_c}(\tfrac{V_{k,t}-k'}{t-k+1}> c). \label{eqn:prop_4_post_nui_bdd}
	\end{align}
	
	From \cref{lem:truncation_error_bound}, there exists $N$ such that $\P{\nu_n,\nu_c}(\left| \overline{L^{k':V_{k,t}-1}_{>N}}\right|>\tfrac{\epsilon}{4})<\tfrac{\delta}{3}$. Next, by choosing $c=\tfrac{\epsilon}{4N}$, we have $\P{\nu_n,\nu_c}(\left\{\tfrac{V_{k,t}-k'}{t-k+1}\left| \overline{L^{k':V_{k,t}-1}_{\leq N}}\right|>\tfrac{\epsilon}{4}\right\}\bigcap\left\{\tfrac{V_{k,t}-k'}{t-k+1} \leq c\right\})=0.$
	Finally, from \cref{lem:no_nui_change_gap_small}, there exists $T_2$ such that for all $t\geq T_2$, we have
	$
	\P{\infty,\nu_c}(\tfrac{V_{k,t}-k'}{t-k+1}> c)<\tfrac{\delta}{3}.$
	The right-hand side of \cref{eqn:prop_4_post_nui_bdd} is then upper bounded by $\delta$, and the proof is complete.

	\section{Proof of Lemma~\ref{lem:fa_prob}}\label[appendix]{sec:AppLem1}
	For any $b>0$, we have
	\begin{align*}
	& \P{\nu_n,\infty}(\eta^1<\infty)                                                                                                                                                                           \\
	& =\sum_{k=1}^\infty\P{\nu_n,\infty}(\eta^1=k)                                                                                                                                                              \\
	& =\sum_{k=1}^\infty\int_{\{\eta^1=k\}}\prod_{i=1}^k h_{\nu_n,\infty,i}(x_i)\  \ud \mathbf{x}_{1:k}                                                                                                         \\
	& \leq \sum_{k=1}^\infty\int_{\{\eta^1=k\}}e^{-b}\Lambda(1,k)\prod_{i=1}^k h_{\nu_n,\infty,i}(x_i)\  \ud \mathbf{x}_{1:k}                                                                                   \\
	\begin{split}& = e^{-b}\sum_{k=1}^\infty\int_{\{\eta^1=k\}}\tfrac{\prod_{i=1}^k g(x_i)}{\max_{1\leq j\leq k+1}\prod_{i=1}^{j-1}f(x_i)\prod_{i=j}^{k}f_n(x_i)}\\&\quad\quad\quad\quad\quad\times\prod_{i=1}^k h_{\nu_n,\infty,i}(x_i)\  \ud \mathbf{x}_{1:k} \end{split}\\
	& \leq e^{-b}\sum_{k=1}^\infty\int_{\{\eta^1=k\}}\prod_{i=1}^k g(x_i)\  \ud \mathbf{x}_{1:k}                                                                                                                \\
	& \leq e^{-b}\P{\infty,1}(\eta^1<\infty)\leq e^{-b}.
	\end{align*}
	The proof that $\P{\nu_n,\infty}(\eta^1_n<\infty)\leq e^{-b}$ is similar. We then have
	$
	\P{\nu_n,\infty}(\min\{\eta^1,\eta^1_n\}<\infty)\leq 2e^{-b}.
	$
	Since $\tauSGLR=\min\{\tau(b),\tau_n(b)\}=\inf_{k\geq1}\min\{\eta^k,\eta^k_n\}$, applying \cite[Theorem~6.16]{poor2009quickest}, we obtain $\E{\nu_n,\infty}[\tauSGLR(b)]\geq \tfrac{1}{2}e^b$ and \cref{ARL_tauWSGLR} follows from $\tauWSGLR\geq\tauSGLR$. The proof is now complete.

	\section{Proof of Lemma~\ref{lem:limsup_assumption}}\label[appendix]{sec:AppLem2}
	It suffices to show that for any $\epsilon>0$, there exists $T$ such that for all $t\geq T$ we have
	\begin{align}
	\sup_{\mathclap{\nu_n\in\mathbb{N},1\leq \nu_c \leq k}}\P{\nu_n,\nu_c}(\tfrac{1}{t}\log \Lambda(k,k+t-1)-I\leq -\delta)\leq\epsilon.\label{eqn:limsup_assumption_epsilon_delta}
	\end{align}
	The set over which the supremum in \cref{eqn:limsup_assumption_epsilon_delta} is taken can be divided into two subsets: $A_1=\{(\nu_n,\nu_c,k)\ :\ \max\{\nu_c,\nu_n\}\leq k\leq t\}$ and $A_2=\{(\nu_n,\nu_c,k)\ :\ \nu_c\leq k<\nu_n\leq t\}$. We have
	\begin{align}
	\begin{split}&\sup_{A_1}\P{\nu_n,\nu_c}(\tfrac{1}{t}\log \Lambda(k,k+t-1)-I\leq -\delta)\\
	& =\P{1,1}(\tfrac{1}{t}\log \Lambda(1,t)-I\leq -\delta)
	\leq \tfrac{\epsilon}{2},\end{split}\label{eqn:reduction_of_supremum_1}
	\end{align}
	where the last inequality follows from \cref{thm:convergence_in_probability_general} for $t$ sufficiently large.
	
	If $\nu_c\leq k<\nu_n\leq t$, we obtain
	\begin{align}
	& \P{\nu_n,\nu_c}(\tfrac{1}{t}\log \Lambda(k,k+t-1)-I\leq -\delta)\nn
	& =\P{\nu_n-k+1,1}(\tfrac{1}{t}\log \Lambda(1,t)-I\leq -\delta) \nn
	& \leq \P{\nu_n-k+1,1}(\tfrac{\nu_n-k}{t}\left(\tfrac{\log \Lambda(1,\nu_n-k)}{\nu_n-k}-I\right)\leq -\tfrac{\delta}{2})\nonumber                  \\
	& \quad+\P{\nu_n-k+1,1}(\tfrac{t-(\nu_n-k)}{t}\left(\tfrac{\log \Lambda(\nu_n-k+1,t)}{t-(\nu_n-k)}-I\right)\leq -\tfrac{\delta}{2})\nonumber       \\
	& \leq \P{\infty,1}(\tfrac{\nu_n-k}{t}\left(\tfrac{\log \Lambda(1,\nu_n-k)}{\nu_n-k}-I\right)\leq -\tfrac{\delta}{2})\label{eqn:case2_infty1}      \\
	& \quad+\P{1,1}(\tfrac{t-(\nu_n-k)}{t}\left(\tfrac{\log \Lambda(1,t-(\nu_n-k))}{t-(\nu_n-k)}-I\right)\leq -\tfrac{\delta}{2}).\label{eqn:case2_11}
	\end{align}
	From \cref{thm:convergence_in_probability_general}, there exists $N_1$ such that for all $n\geq N_1$, we have
	\begin{align}
	\P{\infty,1}(\tfrac{1}{n}\log \Lambda(1,n)-I\leq -\tfrac{\delta}{2})\leq\tfrac{\epsilon}{2},\label{eqn:conv_in_prop_bound_infty1} \\
	\P{1,1}(\tfrac{1}{n}\log \Lambda(1,n)-I\leq -\tfrac{\delta}{2})\leq\tfrac{\epsilon}{2}.\label{eqn:conv_in_prop_bound_11}
	\end{align}
	From Markov's inequality and \cref{assumpt:moments}, there exists $N_2$ such that for all $1\leq n< N_1$ and $m\geq N_2$, we have
	\begin{align}
	\begin{split}&\P{\infty,1}(\tfrac{n}{m}\left(\tfrac{1}{n}\log \Lambda(1,n)-I\right)\leq -\tfrac{\delta}{2})\\
	&\quad\quad\quad\quad \leq\tfrac{2n}{m\delta}\E{\infty,1}[\left|\tfrac{1}{n}\log \Lambda(1,n)-I\right|]
	\leq \tfrac{\epsilon}{2},\end{split}\label{eqn:markov_bound_infty1}                            \\
	\begin{split}&\P{1,1}(\tfrac{n}{m}\left(\tfrac{1}{n}\log \Lambda(1,n)-I\right)\leq -\tfrac{\delta}{2})\\
	&\quad\quad\quad\quad \leq\tfrac{2n}{m\delta}\E{1,1}[\left|\tfrac{1}{n}\log \Lambda(1,n)-I\right|]
	\leq \tfrac{\epsilon}{2}.\end{split}\label{eqn:markov_bound_11}
	\end{align}
	
	Next, we show that for any $t\geq T=2N_1+N_2$, both \cref{eqn:case2_infty1} and \cref{eqn:case2_11} are bounded by $\epsilon/2$. There are three possible cases:
	\begin{enumerate}
		\item$\nu_n-k\geq N_1$ and $t-(\nu_n-k)\geq N_1$,
		\item$\nu_n-k< N_1$ and $t-(\nu_n-k)\geq N_1$,
		\item$\nu_n-k\geq N_1$ and $t-(\nu_n-k)< N_1$.
	\end{enumerate}
	Applying \cref{eqn:conv_in_prop_bound_infty1} and \cref{eqn:conv_in_prop_bound_11} in the first case, \cref{eqn:conv_in_prop_bound_11} and \cref{eqn:markov_bound_infty1} in the second case, and \cref{eqn:conv_in_prop_bound_infty1} and \cref{eqn:markov_bound_11} in the thrid case to \cref{eqn:case2_infty1} and \cref{eqn:case2_11}, respectively, we obtain
	$
	\sup_{A_2} \P{\nu_n,\nu_c}(\tfrac{1}{t}\log \Lambda(k,k+t-1)-I\leq -\delta)\leq\epsilon.
	$
	The proof for \ref{eqn:limsup:it2} is similar and  proof is now complete.

	\section{Proof of Proposition~\ref{prop:ADD}}\label[appendix]{sec:AppProp3}
	From \cref{mb}, there exists $\gamma>0$ such that $m_b \geq (1+\gamma) b/I$ for all $b$ sufficiently large. For any $0<\epsilon < \gamma/(1+\gamma)$, let $n_b=\ceil{\tfrac{b}{(1-\epsilon)I}}$ and $\delta=\epsilon I$. There exists $b_1>0$ such that $n_b(I-\delta) \geq b$ for all $b\geq b_1$. From \cref{lem:limsup_assumption}, by choosing $b_1$ sufficiently large, we have for all $b\geq b_1$,
	\begin{align}
	& \sup_{1\leq \nu_c\leq k}\P{\nu_n,\nu_c}(\log\Lambda(k,k+n_b-1)<b)\nn
	& \leq \sup_{\mathclap{1\leq \nu_c\leq k}}\P{\nu_n,\nu_c}(\log\Lambda(k,k+n_b-1)\leq n_b(I-\delta))\leq \epsilon.\label{eqn:boundnb}
	\end{align}
	Let $b_2 \geq b_1$ be such that $I/b_2 \leq 1+\gamma - (1-\epsilon)^{-1}$. Then, for $b\geq b_2$, we have
	\begin{align*}
	\tfrac{m_b}{n_b}
	& \geq \tfrac{bI^{-1}(1+\gamma)}{\tfrac{bI^{-1}}{1-\epsilon}+1} \\
	& \geq \tfrac{1+\gamma}{(1-\epsilon)^{-1}+I/b}\geq 1.
	\end{align*}
	
	For any $k,\nu_n,\nu_c\geq 1$, we then have
	\begin{align*}
	& \esssup \P{\nu_n,\nu_c}(\widetilde{\tau}_n(b)-\nu_c+1 > kn_b){X_1,\ldots,X_{\nu_c-1}}                                                      \\
	\begin{split}& = \esssup \mathbb{P}_{\nu_n,\nu_c}\left(\max_{t-m_{b}\leq k'\leq t}\log\Lambda(k',t)<b\right.\\& \quad\quad\quad\quad\quad\quad\quad\text{for all $t\leq \nu_c+kn_b-1$}\ \bigg|
	X_1,\ldots,X_{\nu_c-1}\bigg)\end{split} \\
	\begin{split}& \leq \esssup \mathbb{P}_{\nu_n,\nu_c}(\log\Lambda(\nu_c+(j-1)n_b,\nu_c+jn_b-1)<b\\& \quad\quad\quad\quad\quad\quad\quad\quad\quad\quad\text{for all $1\leq j \leq k$}\ |\ X_1,\ldots,X_{\nu_c-1})    \end{split}  \\
	& =\prod_{j=1}^k \P{\nu_n,\nu_c}(\log\Lambda(\nu_c+(j-1)n_b,\nu_c+jn_b-1)<b)\leq \epsilon^k,
	\end{align*}
	where the last equality follows from independence and the last inequality from \cref{eqn:boundnb}. Therefore, for any $b\geq b_2$, we have
	\begin{align*}
	& \sup_{{\nu_n,\nu_c\geq 1}}\esssup \E{\nu_n,\nu_c}[(\widetilde{\tau}_n(b)-\nu_c+1)^+]{X_1,\ldots,X_{\nu_c-1}}                                            \\
	& =\sup_{{\nu_n,\nu_c\geq 1}}\esssup \sum_{i=0}^\infty\P{\nu_n,\nu_c}((\widetilde{\tau}_n(b)-\nu_c+1)>i){X_1,\ldots,X_{\nu_c-1}}                          \\
	& \leq\sum_{k=0}^\infty n_b \sup_{{\nu_n,\nu_c\geq 1}}\esssup \P{\nu_n,\nu_c}(\widetilde{\tau}_n(b)-\nu_c+1>kn_b){X_1,\ldots,X_{\nu_c-1}}                 \\
	& \leq \tfrac{n_b}{1-\epsilon}
	\leq \tfrac{b}{I(1-\epsilon)^2}+\tfrac{1}{1-\epsilon}=b\left(I^{-1}+\tfrac{2\epsilon-\epsilon^2}{I(1-\epsilon)^2}+\tfrac{1}{b(1-\epsilon)}\right),
	\end{align*}
	which yields \ref{prop:ADD:it1}. The proof for \ref{prop:ADD:it2} is similar and the proposition is proved.

	\section{Proof of Theorem~\ref{thm:main_result}}\label[appendix]{sec:AppThm1}
	
	From \cref{lem:fa_prob}, taking infimum on both sides of \cref{ARL_tauWSGLR}, we obtain $
	\ARL(\tau_{\text{W-SGLR}}(b))=\inf_{\nu_n\in\mathbb{N}\cup\{\infty\}}\E{\nu_n,\infty}[\tau_{\text{W-SGLR}}(b)]\geq \tfrac{1}{2}e^b.
	$
	Since $\tau_{\text{W-SGLR}}\leq\widetilde{\tau}$ and $\tau_{\text{W-SGLR}}\leq\widetilde{\tau}_n$, by \cref{prop:ADD}, we have
	$
	\text{WADD}(\tau_{\text{W-SGLR}}(b))\leq(I^{-1}+o(1))b
	$
	as $b\to\infty$.
	
	To see that $\tau_{\text{W-SGLR}}(b)$ is asymptotically optimal when \cref{assumpt:kldiv} is satisfied, let $C_{\gamma}=\{\tau\ :\ \ARL(\tau)\geq \gamma\}$ be the set of stopping times satisfying $\ARL(\tau)\geq \gamma$. By expanding $\WADD(\tau)$ using \cref{eqn:WADD}, we obtain
	\begin{align}
	\inf_{\tau\in C_{\gamma}}\WADD(\tau)&=\inf_{\tau\in C_{\gamma}}\sup_{\nu_n\in\mathbb{N}\cup\{\infty\}}\WADD_{\nu_n}(\tau)\nn
	&\geq \sup_{\nu_n\in\mathbb{N}\cup\{\infty\}}\inf_{\tau\in C_{\gamma}}\WADD_{\nu_n}(\tau)\label{eqn:min-max_ineq}\\
	&\geq \sup_{\nu_n\in\{0,\nu_c,\infty\}}\inf_{\tau\in C_{\gamma}}\WADD_{\nu_n}(\tau),\label{eqn:supinf}
	\end{align} 
	where \cref{eqn:min-max_ineq} is due to the min-max inequality\cite{boyd2004convex}.  For each of the cases $\nu_n\in\{0,\nu_c,\infty\}$, by Theorem 6.17 in \cite{poor2009quickest}, we have 
	\begin{align*}
	\inf_{\tau\in C_{\gamma}}\WADD_{\nu_n}(\tau) & \geq\left(\KLD{g_n}{f_n}^{-1}+o(1)\right)b, \quad \text{when $\nu_n=0$,}   \\
	\inf_{\tau\in C_{\gamma}}\WADD_{\nu_n}(\tau) & \geq\left(\KLD{g_n}{f}^{-1}+o(1)\right)b, \quad \text{when $\nu_n=\nu_c$,} \\
	\inf_{\tau\in C_{\gamma}}\WADD_{\nu_n}(\tau) & \geq\left(\KLD{g}{f}^{-1}+o(1)\right)b, \quad \text{when $\nu_n=\infty$.}
	\end{align*}
	Since \cref{assumpt:kldiv} is satisfied, we have 
	\begin{align}
	I=\min\left\{\KLD{g}{f},\ \KLD{g_n}{f},\ \KLD{g_n}{f_n}\right\}.\label{eqn:I}
	\end{align} 
	Therefore, from \cref{eqn:supinf} and \cref{eqn:I}, we obtain
	\begin{align*}
	\inf_{\tau\in C_{\gamma}}\WADD_{\nu_n}(\tau)
	& \geq \Big(\max\big\{\KLD{g_n}{f_n}^{-1},\ \KLD{g_n}{f}^{-1}, \\
	& \qquad\qquad \KLD{g}{f}^{-1}\big\} + o(1) \Big)b \\
	& = \left(I^{-1}+o(1)\right)b,
	\end{align*}
	and the proof is now complete.

	\section{Proof of Lemma~\ref{lem:error_prob_glrt}}\label{sec:App_Lemma_error_prob}

	We use techniques is similar to \cite{lai98} to prove Lemma~\ref{lem:error_prob_glrt}. To analyze the probability $\P{\nu_n,\infty}(\widehat{\eta}_k<\infty)$, we use a change-of-measure argument. For any $0<\delta<1$, choose $b_\delta\geq0$ so that for any $b\geq b_\delta$
	\begin{align*}|\Theta| \pi^{-\tfrac{d}{2}}{\Gamma\left(\tfrac{d}{2}+1\right)}b^{\tfrac{d}{2}} \exp\left(-b\right) \leq  \exp\left(-(1-\delta)b\right),
	\end{align*}
	where $|\Theta|$ is the volume or Lebesgue measure of $\Theta\subset\Real^d$ and $\Gamma(\cdot)$ is the gamma function. From Kolmogorov's Consistency Theorem, there is a probability measure $G_\theta$ for the stochastic process $(X_i)_{i\geq k}$ under which the pdf of each $X_i$ is $g(\cdot;\theta)$.  Define a measure $H(\cdot)=\int_{\Theta}G_\theta(\cdot)\ \ud \theta$. Since $\Theta$ is compact in $\mathbb{R}^d$, the measure $H$ is finite. For each $t\geq k$, the Radon-Nikodym derivative of the law of $(X_k,X_{k+1},\ldots,X_t)$ under $H$ \gls{wrt} $\P{\nu_n,\infty}$ is
	\begin{align*}
	R_t=\int_{\Theta}\exp\left(\sum_{i=k}^t\log \tfrac{g(X_i;\theta)}{h_{\nu_n,\infty,\theta,\theta_n,i}(X_i)}\right)\ \ud\theta,
	\end{align*}
	which follows from Fubini's Theorem. By Wald's likelihood ratio identity,
	\begin{align}
	\P{\nu_n\infty}(\widehat{\eta}_k<\infty) & =\int_{\{\widehat{\eta}_k<\infty\}}R_{\widehat{\eta}_k}^{-1}\ \ud H\nn
	& =\int_{\Theta}\left\{\int_{\{\widehat{\eta}_k<\infty\}}R_{\widehat{\eta}_k}^{-1}\ \ud G_\theta\right\}\ \ud \theta. \label{eqn:wald_llr_identity}
	\end{align}
	Suppose $\widehat{\eta}_k=t$. Since $\widehat{\theta}=\argmax_\theta \log\widehat{\Lambda}(k,t,\theta)\in\text{Int}(\Theta)$, from Taylor series, there exists $\theta^*\in\Theta$ such that
	\begin{align*}
	\sum_{i=k}^t \log g(x_i;\theta)=\sum_{i=k}^t \log g(x_i;\hat{\theta})+\tfrac{1}{2}(\theta-\hat{\theta})^T\left[\nabla^2\sum_{i=k}^t \log g(x_i,\theta^*)\right](\theta-\hat{\theta}).
	\end{align*}
	Thus, for $\|\theta-\widehat{\theta}\|<1/\sqrt{b}$, we have
	\begin{align*}
	\log\widehat{\Lambda}(k,t,\hat{\theta})-\log\widehat{\Lambda}(k,t,\theta) 
	& =-\tfrac{1}{2}(\theta-\hat{\theta})^T\left[\nabla^2 \log\widehat{\Lambda}(k,t,\theta^*)\right](\theta-\hat{\theta}) \\
	& \leq\tfrac{1}{2}\|\theta-\widehat{\theta}\|^2\lambda_{\max}\left(-\nabla^2\log\widehat{\Lambda}(k,t,\theta)\right)\\
	&\leq \ofrac{2}\|\theta-\widehat{\theta}\|^2 b,
	\end{align*}
	where the last inequality follows from $\sup_{\|\theta-\widehat{\theta}\|<1/\sqrt{b}}\lambda_{\max}\left(-\nabla^2\log\widehat{\Lambda}(k,t,\theta)\right)\leq b$. We obtain
	\begin{align*}
	\tfrac{R_t}{\widehat{\Lambda}(k,t,\widehat{\theta})} & =\int_{\Theta} \exp\left(-\left[\log\widehat{\Lambda}(k,t,\widehat{\theta})-\sum_{i=k}^t \log \tfrac{g(X_i,\theta)}{h_{\nu_n,\infty,\theta,\theta_n,i}(X_i)} \right]\right)\ \ud \theta \\
	& \geq\int_{\Theta} \exp\left(-\left[\log\widehat{\Lambda}(k,t,\widehat{\theta})-\log \widehat{\Lambda}(k,t,\theta) \right]\right)\ \ud \theta                                           \\
	& \geq\int_{\|\theta-\widehat{\theta}\|\leq 1/\sqrt{b}} \exp\left(-\ofrac{2}\|\theta-\widehat{\theta}\|^2b\right)\ \ud \theta\\
	& \geq\int_{\|\theta-\widehat{\theta}\|\leq 1/\sqrt{b}} 1\ \ud \theta \\
	& \geq\tfrac{\pi^{\tfrac{d}{2}}}{\Gamma(\tfrac{d}{2}+1)}b^{-\tfrac{d}{2}} 
	\end{align*}
	
	Therefore, we have
	\begin{align*}
	R_t & \geq \tfrac{\pi^{\tfrac{d}{2}}}{\Gamma(\tfrac{d}{2}+1)}b^{-\tfrac{d}{2}}\widehat{\Lambda}(k,t,\widehat{\theta}) \\
	& \geq \tfrac{\pi^{\tfrac{d}{2}}}{\Gamma(\tfrac{d}{2}+1)}b^{-\tfrac{d}{2}} \exp\left(b\right).
	\end{align*}
	This yields the upper bound 
	\begin{align*}
	R_t^{-1} & \leq \pi^{-\tfrac{d}{2}}{\Gamma(\tfrac{d}{2}+1)}b^{\tfrac{d}{2}} \exp\left(-b\right).
	\end{align*}	
	Applying this upper bound to \eqref{eqn:wald_llr_identity}, we obtain
	\begin{align*}
	\P{\nu_n,\infty}(\widehat{\eta}_k<\infty)
	& \leq  \int_{\Theta}\left\{ \pi^{-\tfrac{d}{2}}{\Gamma\left(\tfrac{d}{2}+1\right)}b^{\tfrac{d}{2}} \exp\left(-b\right)\int_{\{\widehat{\eta}_k<\infty\}}\ \ud G_\theta\right\}\ \ud \theta \\
	& \leq  \int_{\Theta}\left\{ \pi^{-\tfrac{d}{2}}{\Gamma\left(\tfrac{d}{2}+1\right)}b^{\tfrac{d}{2}} \exp\left(-b\right) \right\}\ \ud \theta                                                \\
	& =  |\Theta| \pi^{-\tfrac{d}{2}}{\Gamma\left(\tfrac{d}{2}+1\right)}b^{\tfrac{d}{2}} \exp\left(-b\right) \\
	& \leq  \exp\left(-(1-\delta)b\right),
	\end{align*}
	for all $b\geq b_\delta$. The proof that $\P{\nu_n,\infty}(\widehat{\eta}_{n,l}<\infty)\leq \exp\left(-(1-\delta)b\right)$ is similar, and the lemma is proved.
	
	\bibliographystyle{IEEEtran}
	\bibliography{IEEEabrv,StringDefinitions,refs}
	
\end{document}